 \documentclass[12pt]{article}
 \usepackage{amsfonts}
 \usepackage{mathrsfs}
 \usepackage{mathrsfs,amsfonts}
 \usepackage{color}
 \usepackage{amsfonts,amssymb,amsmath,mathrsfs}
 \usepackage[svgnames]{xcolor}

 \newtheorem{theorem}{Theorem}[section]
 \newtheorem{lemma}[theorem]{Lemma}
 \newtheorem{corol}[theorem]{Corollary}
 \newtheorem{prop}[theorem]{Proposition}
 \newtheorem{remark}[theorem]{Remark}
 \newtheorem{example}[theorem]{Example}
 \newtheorem{condition}[theorem]{Condition}

 \def\btheorem{\begin{theorem}\sl{}
 \def\etheorem{\end{theorem}}}
 \def\blemma{\begin{lemma}\sl{}
 \def\elemma{\end{lemma}}}
 \def\bcorollary{\begin{corol}\sl{}
 \def\ecorollary{\end{corol}}}
 \def\bproposition{\begin{prop}\sl{}
 \def\eproposition{\end{prop}}}

 \def\benumerate{\begin{enumerate}}\def\eenumerate{\end{enumerate}}
 \def\bitemize{\begin{itemize}}\def\eitemize{\end{itemize}}
 \def\itm{\item}

 \def\beqlb{\begin{eqnarray}}
 \def\eeqlb{\end{eqnarray}}
 \def\beqnn{\begin{eqnarray*}}
 \def\eeqnn{\end{eqnarray*}}

 \def\qqquad{\qquad\qquad}

 \def\proof{\noindent{\it Proof.~~}}
 \def\qed{\hfill$\Box$\medskip}

 \def\<{\langle}\def\>{\rangle}

 \def\mcr{\mathscr}\def\mbb{\mathbb}
 \def\mbf{\mathbf}\def\mrm{\mathrm}

 \def\ar{\!\!&}\def\nnm{\nonumber}\def\ccr{\nnm\\}

 \def\d{\mrm{d}}\def\e{\mrm{e}}
 \def\const{\mathrm{const}}

\begin{document}

$~$

\bigskip\bigskip

\centerline{\Large\bf Construction of continuous-state branching}

\smallskip

\centerline{\Large\bf processes in varying environments\,\footnote{Supported by the National Natural Science Foundation of China (No.11531001).}}

\bigskip

\centerline{Rongjuan Fang and Zenghu Li}

\bigskip

\centerline{School of Mathematical Sciences, Beijing Normal University,}

\centerline{Beijing 100875, People's Republic of China}

\centerline{E-mails: {\tt fangrj@mail.bnu.edu.cn} and {\tt lizh@bnu.edu.cn}}

\bigskip\bigskip

{\narrower{\narrower

\noindent\textit{Abstract}: A continuous-state branching process in varying environments is constructed by the pathwise unique solution to a stochastic integral equation driven by time-space noises. The process arises naturally in the limit theorem of Galton--Watson processes in varying environments established by Bansaye and Simatos (2015). In terms of the stochastic equation we clarify the behavior of the continuous-state process at its bottlenecks, which are the times when it arrives at zero almost surely by negative jumps.

\medskip



\par}\par}

\bigskip\bigskip


\section{Introduction}

\setcounter{equation}{0}

\textit{Continuous-state branching processes} (CB-processes) are often used to model the stochastic evolution of large populations with small individuals. The \textit{branching property} means intuitively that different individuals in the population propagate independently of each other. The study of such processes was initiated by Feller (1951), who noticed that a diffusion process may arise in a limit theorem of rescaled \textit{Galton--Watson branching processes} (GW-processes). The basic structures of general CB-processes were discussed in Ji\v{r}ina (1958). It was proved in Lamperti (1967a) that the class of CB-processes with homogeneous transition semigroups coincides with that of scaling limits of classical GW-processes; see also Aliev and Shchurenkov (1982) and Grimvall (1974). The connection between CB-processes and time changed L\'{e}vy processes was established by Lamperti (1967b). A general existence theorem for homogeneous CB-processes was proved in Silverstein (1968); see also Watanabe (1969) and Rhyzhov and Skorokhod (1970). The approach of stochastic equations for CB-processes without or with immigration has been developed by Bertoin and Le~Gall (2006), Dawson and Li (2006, 2012), Fittipaldi and Fontbona (2012), Fu and Li (2010), Li (2011, 2019+), Pardoux (2016) and many others.

There have also been some attempts at the understanding of inhomogeneous CB-processes. Let $X= \{X(t): t\in I\}$ be a Markov process with state space $[0,\infty]$ and inhomogeneous transition semigroup $\{Q_{r,t}: t\ge r\in I\}$, where $I\subset \mbb{R}$ is an interval. We call $X$ a \textit{CB-process in varying environments} (CBVE-process) if there is a family of continuous mappings $\{v_{r,t}: t\ge r\in I\}$ on $(0,\infty)$ so that
 \beqlb\label{eq1.1}
\int_{[0,\infty]}\e^{-\lambda y}Q_{r,t}(x,\d y)= \e^{-xv_{r,t}(\lambda)}, \qquad \lambda> 0, x\in [0,\infty]
 \eeqlb
with $\e^{-\lambda y}= 0$ for $y=\infty$ by convention. It is natural to expect that the processes defined by \eqref{eq1.1} are scaling limits of \textit{GW-processes in varying environments} (GWVE-processes), where individuals in different generations may have different reproduction distributions. The understanding of the CBVE-processes is important since they provide the bases of further study of \textit{CB-processes in random environments} (CBRE-processes). The reader may refer to Bansaye et al.\ (2013, 2019), Bansaye and Simatos (2015), He et al.\ (2018), Helland (1981), Kurtz (1978), Li and Xu (2018), Palau et al.\ (2016), Palau and Pardo (2017, 2018) and the references therein for some progresses in the study. In particular, a scaling limit theorem for a sequence of GWVE-processes was proved by Bansaye and Simatos (2015), who provided a general sufficient condition for the weak convergence of the sequence and showed a CBVE-process indeed arises as the limit. Their condition allows infinite variance of the reproduction distributions and extends considerably the results in this line established before. But the general existence theorem for the CBVE-process was not provided in Bansaye and Simatos (2015). In fact, with their approach they need to avoid the \textit{bottlenecks}, which are the times when the process arrives at zero a.s.\ by negative jumps. The determination of the behavior of the CBVE-process at the bottlenecks was left open in Bansaye and Simatos (2015).

The purpose of this work is to give a construction of the CBVE-process under reasonably general assumptions and clarify its behavior at the bottlenecks. Let $b_1$ and $c$ be c\`{a}dl\`{a}g functions on $[0,\infty)$ satisfying $b_1(0)= c(0)= 0$ and having locally bounded variations. Let $m$ be a $\sigma$-finite measure on $(0,\infty)^2$ satisfying
 \beqlb\label{eq1.2}
m_1(t):= \int_0^t\int_0^\infty (1\wedge z^2) m(\d s,\d z)< \infty, \qquad t\ge 0.
 \eeqlb
Here and in the sequel, we understand, for $t\ge r\in \mbb{R}$,
 \beqnn
\int_r^t = \int_{(r,t]} = -\int_{t}^{r},
 \qquad
\int_r^\infty = \int_{(r,\infty)} = -\int_{\infty}^{r}.
 \eeqnn
Let us consider the backward integral evolution equation:
 \beqlb\label{eq1.3}
v_{r,t}(\lambda)= \lambda - \int_r^t v_{s,t}(\lambda)b_1(\d s) - \int_r^t v_{s,t}(\lambda)^2 c(\d s) - \int_r^t\int_0^\infty K_1(v_{s,t}(\lambda),z) m(\d s,\d z),
 \eeqlb
where $K_1(\lambda,z)= \e^{-\lambda z} - 1 + \lambda z1_{\{z\le 1\}}$. This is an equivalent reformulation of the equation in Theorem~2.2 of Bansaye and Simatos (2015). We say the parameters $(b_1,c,m)$ are \textit{weakly admissible} provided:
 \bitemize

\itm[(1.A)] $t\mapsto c(t)$ is increasing and continuous;

\itm[(1.B)] for every $t>0$ we have
 \beqlb\label{eq1.4}
\Delta b_1(t) + \int_0^1 z m(\{t\},\d z)\le 1,
 \eeqlb
where $\Delta b_1(t)= b_1(t)-b_1(t-)$.

 \eitemize
It is natural to introduce condition \eqref{eq1.4} to ensure that the solution of \eqref{eq1.3} stays positive ($=$\,nonnegative). In fact, from \eqref{eq1.3} we have, for any $0<\varepsilon\le 1$,
 \beqnn
v_{t-,t}(\lambda)\ar=\ar \lambda[1-\Delta b_1(t)] - \int_0^\infty K_1(\lambda,z) m(\{t\},\d z) \cr
 \ar\le\ar
\lambda\bigg[1 - \Delta b_1(t) - \int_\varepsilon^1 z m(\{t\},\d z)\bigg] + m(\{t\}\times (\varepsilon,\infty)),
 \eeqnn
and hence $v_{t-,t}(\lambda)< 0$ for sufficiently large $\lambda> 0$ if \eqref{eq1.4} is not satisfied. Let $J= \{s>0: \Delta b_1(s)= 1\}$ and $K= \{s\in J: m(\{s\}\times (0,\infty))= 0\}$. We say the parameters $(b_1,c,m)$ are \textit{admissible} if they are weakly admissible and $K$ is an empty set.

\btheorem\label{th1.1} Let $(b_1,c,m)$ be admissible parameters. Then for $t\ge 0$ and $\lambda> 0$ there is a unique bounded and strictly positive solution $[0,t]\ni r\mapsto v_{r,t}(\lambda)$ to the integral evolution equation \eqref{eq1.3} and a transition semigroup $(Q_{r,t})_{t\ge r}$ on $[0,\infty]$ is defined by \eqref{eq1.1}. \etheorem

\btheorem\label{th1.2} Let $(b_1,c,m)$ be admissible parameters. Then for any $t\ge 0$, $r\mapsto v_{r,t}(0):= \lim_{\lambda\downarrow 0} v_{r,t}(\lambda)$ is the largest positive solution to \eqref{eq1.3} with $\lambda= 0$ and $r\mapsto v_{r,t}(\infty):= \lim_{\lambda\uparrow \infty} v_{r,t}(\lambda)$ is the smallest positive solution to \eqref{eq1.3} with $\lambda= \infty$. \etheorem

Suppose that $(\Omega, \mcr{F}, \mcr{F}_t, \mbf{P})$ is a filtered probability space satisfying the usual hypotheses. Let $W(\d s,\d u)$ be a time-space $(\mcr{F}_t)$-Gaussian white noise on $(0,\infty)^2$ with intensity $2c(\d s)\d u$. Let $M(\d s,\d z,\d u)$ be a time-space $(\mcr{F}_t)$-Poisson random measure on $(0,\infty)^3$ with intensity $m(\d s,\d z)\d u$. Denote by $\tilde{M}(\d s,\d z,\d u)$ the compensated measure of $M(\d s,\d z,\d u)$. Given an $\mcr{F}_0$-measurable random variable $X(0)\ge 0$, we consider the stochastic integral equation:
 \beqlb\label{eq1.5}
X(t)\ar=\ar X(0) + \int_0^t\int_0^{X(s-)}W(\d s,\d u) + \int_0^t\int_0^1 \int_0^{X(s-)} z \tilde{M}(\d s,\d z,\d u), \cr
 \ar\ar\quad
- \int_0^tX(s-)b_1(\d s) + \int_0^t\int_1^\infty\int_0^{X(s-)}z M(\d s,\d z,\d u).
 \eeqlb
By saying the positive c\`{a}dl\`{a}g process $\{X(t): t\ge 0\}$ in $[0,\infty]$ is a \textit{solution} to \eqref{eq1.5} we mean the equation holds a.s.\ if $t$ is replaced by $t\land\tau_k$ for every $t\ge 0$, where $\tau_k= \inf\{t\ge 0: X(t)\ge k\}$, and the states $0$ and $\infty$ are traps for $\{X(t): t\ge 0\}$.

\btheorem\label{th1.3} Let $(b_1,c,m)$ be admissible parameters. Then there is a pathwise unique solution $\{X(t): t\ge 0\}$ to \eqref{eq1.5} and the solution is a CBVE-process with transition semigroup $(Q_{r,t})_{t\ge r}$ defined by \eqref{eq1.1} and \eqref{eq1.3}.
\etheorem

The CBVE-process constructed by \eqref{eq1.3} and \eqref{eq1.5} is a generalization of the model studied in Ji\v{r}ina (1958), where a smoothness was assumed for \eqref{eq1.3}. We shall first treat special forms of \eqref{eq1.3} and \eqref{eq1.5} by imposing an integrability condition stronger than \eqref{eq1.2}, which implies the CBVE-process has finite first moments. The existence of the cumulant semigroup is constructed by an iteration argument combined with an inhomogeneous nonlinear $h$-transformation. A suitably chosen transformation of this type changes the CBVE-process into a positive martingale and plays an important role in the establishment of the stochastic equation under the first moment assumption. The solutions to the general equations \eqref{eq1.3} and \eqref{eq1.5} are then obtained by increasing limits. The Poisson random measure in \eqref{eq1.5} does not fit immediately into the framework of single valued point processes developed in standard references such as Ikeda and Watanabe (1989), Jacod and Shiryaev (2003) and Situ (2005). In fact, at a fixed discontinuity $t>0$ the jump size $\Delta X(t)$ of the CBVE-process is identified by a composite L\'{e}vy--It\^{o} representation as the position at \textit{time} $X(t-)$ of a spectrally positive L\'{e}vy process constructed from the random measure $M(\{t\},\d z,\d u)$, which typically has infinitely many atoms. This is essentially different from its homogeneous version discussed in Bertoin and Le~Gall (2006) and Dawson and Li (2006, 2012), where $M(\{t\},\d z,\d u)$ has no more than one atom. The complexity of jumps of the solution makes the treatment of \eqref{eq1.5} much more difficult than the homogeneous equations. The time-space noises in the stochastic equation yield natural interactions among the solutions started from different initial states, which are essential in the analysis of the model. By Theorem~\ref{th1.2}, the uniqueness of solutions to \eqref{eq1.3} holds for $\lambda\ge 0$ if and only if it holds for $\lambda= 0$. This verifies an observation of Rhyzhov and Skorokhod (1970, p.706) in our setting. The probabilistic meanings of the quantities $v_{r,t}(0)$ and $v_{r,t}(\infty)$ are given in \eqref{eq2.16} and \eqref{eq2.17}, respectively.

Let $(b_1,c,m)$ be weakly admissible parameters. We call any moment $s\in K$ a \textit{bottleneck} following the terminology of Bansaye and Simatos (2015). Since $b_1$ is a c\`{a}dl\`{a}g function, we can rearrange $K$ into an increasing (finite or infinite) sequence $\{s_1,s_2, \cdots\}$. For $t> 0$ let $\wp(t)= \max\{s\in K: s\le t\}$ with $\max\emptyset =0$ by convention. By Theorem~\ref{th1.1}, for any $\lambda> 0$ there is a unique bounded and strictly positive solution $r\mapsto v_{r,t}(\lambda)$ to \eqref{eq1.3} on the interval $[\wp(t),t]$. By setting $v_{r,t}(\lambda)= 0$ for $0\le r< \wp(t)$ we can extend $r\mapsto v_{r,t}(\lambda)$ into a solution to \eqref{eq1.3} on $[0,t]$. In this case, we may not be able to define the whole transition semigroup $\{Q_{r,t}: t\ge r\in [0,\infty)\}$ simultaneously by \eqref{eq1.1}. However, for each $i=0,1,2,\cdots$ we can use \eqref{eq1.1} and \eqref{eq1.3} to define a transition semigroup $\{Q_{r,t}: t\ge r\in [s_i,s_{i+1})\}$ on $[0,\infty]$, where we understand $s_0=0$.

In terms of the stochastic equation, the behavior of the CBVE-process at the bottlenecks is clarified as follows. For weakly admissible parameters, we can use Theorem~\ref{th1.3} to see there is still a pathwise unique solution $\{X_0(t): t\ge 0\}$ to \eqref{eq1.5} and its restriction to the time interval $[0,s_1)$ is a CBVE-process with transition semigroup $\{Q_{r,t}: t\ge r\in [0,s_1)\}$. Let $\tau_{0,k}= \inf\{t\ge 0: X_0(t)\ge k\}$ and let $\tau_{0,\infty}= \lim_{k\to \infty}\tau_{0,k}$ be the \textit{explosion time} of $\{X_0(t): t\ge 0\}$. Then we have $X_0(s_1)= 0$ on the event $\{s_1< \tau_{0,\infty}\}$ and $X_0(s_1)= \infty$ on the event $\{\tau_{0,\infty}\le s_1\}$. In fact, for any $r_i\in [s_i,s_{i+1})$, $i=1,2,\cdots$, given the initial value $X_i(r_i)\ge 0$, we can construct a process $\{X_i(t): t\ge r_i\}$ by the pathwise unique solution to a time-shift of \eqref{eq1.5}. The restriction of the solution to $[r_i,s_{i+1})$ is a CBVE-process with transition semigroup $\{Q_{r,t}: t\ge r\in [r_i,s_{i+1})\}$. The behavior of $\{X_i(t): t\in [r_i,s_{i+1})\}$ at $s_{i+1}\in K$ is similar to that of $\{X_0(t): t\in [0,s_1)\}$ at $s_1\in K$.

The remaining part of the paper is organized as follows. In Section~2, some preliminary results are presented. In Section~3, we exploit the existence and uniqueness of solutions to some special cases of \eqref{eq1.3}. The corresponding CBVE-process is constructed in Section~4 by solving a special form of \eqref{eq1.5}. The general results for admissible parameters are proved in Section~5.

\bigskip

\noindent\textbf{Acknowledgements~} We would like to thank Peisen Li for helpful comments on an earlier version of the work. We are grateful to the Laboratory of Mathematics and Complex Systems (Ministry of Education) for providing us the research facilities to carry out the project.


\section{Preliminaries}

\setcounter{equation}{0}

Given a c\`{a}dl\`{a}g function $\alpha$ on $[0,\infty)$ with locally bounded variations and $\alpha(0)= 0$, we write $\Delta\alpha(t)= \alpha(t) - \alpha(t-)$ for the size of its jump at $t>0$ and $\|\alpha\|(t)$ for the total variation of $\alpha$ on $[0,t]$. It is well-known the set $J_\alpha:= \{s> 0: \Delta\alpha(s)\neq 0\}$ is at most countable. The function $\alpha$ has the decomposition $\alpha(t)= \alpha_c(t) + \alpha_d(t)$, where $\alpha_d(t):= \sum_{0<s\le t} \Delta\alpha(s)$ is the \textit{jump part} and $\alpha_c(t):= \alpha(t)-\alpha_d(t)$ is the \textit{continuous part}.


\bproposition\label{th2.2} Suppose that $\alpha$ and $G$ are c\`{a}dl\`{a}g functions on $[0,\infty)$ with locally bounded variations such that $\Delta\alpha(t)> -1$ for every $t> 0$. Let $\zeta$ be the c\`{a}dl\`{a}g function on $[0,\infty)$ such that $\zeta_c(t)= \alpha_c(t)$ and $\Delta\zeta(t)= \log [1+\Delta\alpha(t)]$ for every $t> 0$. Then we have:
 \bitemize

\itm[{\rm(i)}] {\rm(Forward equation)} There is a unique locally bounded solution to:
 \beqlb\label{eq2.6}
F(t)= G(t) + \int_0^t F(s-)\alpha(\d s), \qquad t\ge 0,
 \eeqlb
which is given by
 \beqlb\label{eq2.7}
F(t)= \e^{\zeta(t)-\zeta(0)} G(0) + \int_0^t\e^{\zeta(t)-\zeta(s)}G(\d s).
 \eeqlb

\itm[{\rm(ii)}] {\rm(Backward equation)} For every $t\ge 0$ there is a unique bounded solution to:
 \beqlb\label{eq2.8}
H(r)= G(r) + \int_r^t H(s)\alpha(\d s), \qquad r\in [0,t],
 \eeqlb
which is given by
 \beqlb\label{eq2.9}
H(r)= \e^{\zeta(t)-\zeta(r)} G(t) - \int_r^t\e^{\zeta(s-)-\zeta(r)}G(\d s).
 \eeqlb

 \eitemize
\eproposition

\proof The uniqueness of the solution to \eqref{eq2.6} or \eqref{eq2.8} follows by standard applications of Gronwall's inequalities and is left to the reader. By \eqref{eq2.7} and integration by parts, we have
 \beqnn
F(t)\ar=\ar G(0) + \e^{-\zeta(0)}G(0)\int_0^t \d\e^{\zeta(s)} + \int_0^t G(\d s) + \int_0^t \d\e^{\zeta(s)}\int_0^{s-} \e^{-\zeta(v)}G(\d v) \cr
 \ar=\ar
\e^{-\zeta(0)}G(0)\int_0^t \e^{\zeta(s)}\d\zeta_c(s) + \e^{-\zeta(0)}G(0) \sum_{s\in(0,t]} (\e^{\zeta(s)}-\e^{\zeta(s-)}) + G(t) \cr
 \ar\ar
+ \int_0^t\e^{\zeta(s)} \d\zeta_c(s) \int_0^{s-}\e^{-\zeta(v)}G(\d v) + \sum_{s\in(0,t]} (\e^{\zeta(s)}-\e^{\zeta(s-)}) \int_0^{s-}\e^{-\zeta(v)}G(\d v) \cr
 \ar=\ar
\e^{-\zeta(0)}G(0)\int_0^t \e^{\zeta(s-)}\d\zeta_c(s) + \e^{-\zeta(0)}G(0) \sum_{s\in(0,t]} \e^{\zeta(s-)}(\e^{\Delta\zeta(s)}-1) + G(t) \cr
 \ar\ar
+ \int_0^t\e^{\zeta(s-)} \d\zeta_c(s) \int_0^{s-}\e^{-\zeta(v)}G(\d v) + \sum_{s\in(0,t]} \e^{\zeta(s-)}(\e^{\Delta\zeta(s)}-1) \int_0^{s-}\e^{-\zeta(v)}G(\d v) \cr
 \ar=\ar
\e^{-\zeta(0)}G(0)\int_0^t \e^{\zeta(s-)}\d\alpha(s) + G(t) + \int_0^t\e^{\zeta(s-)} \d\alpha(s) \int_0^{s-}\e^{-\zeta(v)}G(\d v) \cr
 \ar=\ar
G(t) + \int_0^t F(s-) \d\alpha(s).
 \eeqnn
Then $t\mapsto F(t)$ is a solution to \eqref{eq2.6}. Similarly, by \eqref{eq2.9} and integration by parts,
 \beqnn
H(t)\ar=\ar H(r) + \e^{\zeta(t)}G(t) \int_r^t \d\e^{-\zeta(s)} + \int_r^t G(\d s) - \int_r^t\d\e^{-\zeta(s)} \int_s^t\e^{\zeta(v-)}G(\d v) \cr
 \ar=\ar
H(r) - \e^{\zeta(t)}G(t)\int_r^t \e^{-\zeta(s)}\alpha(\d s) + \int_r^t G(\d s) + \int_r^t\e^{-\zeta(s)} \d\alpha(s) \int_s^t\e^{\zeta(v-)}G(\d v) \cr
 \ar=\ar
H(r) - \int_r^t H(s)\alpha(\d s) + G(t) - G(r).
 \eeqnn
Then $r\mapsto H(r)$ is a solution to \eqref{eq2.8} on $[0,t]$. \qed

\bcorollary\label{th2.3} Let $\alpha$ be a c\`{a}dl\`{a}g function on $[0,\infty)$ with locally bounded variations such that $\Delta\alpha(t)> -1$ for $t> 0$. Then for $\lambda\in \mbb{R}$ we have:
 \bitemize

\itm[{\rm(i)}] {\rm(Forward equation)} There is a unique locally bounded solution to:
 \beqlb\label{eq2.10}
\pi_t(\lambda)= \lambda + \int_0^t\pi_{s-}(\lambda)\alpha(\d s), \qquad t\ge 0,
 \eeqlb
which is given by
 \beqlb\label{eq2.11}
\pi_t(\lambda)= \lambda\prod_{s\in(0,t]} (1+\Delta\alpha(s)) \exp\{\alpha_c(t)-\alpha_c(0)\}.
 \eeqlb

\itm[{\rm(ii)}] {\rm(Backward equation)} For every $t\ge 0$ there is a unique bounded solution to:
 \beqlb\label{eq2.12}
\pi_{r,t}(\lambda)= \lambda + \int_r^t\pi_{s,t}(\lambda)\alpha(\d s), \qquad r\in [0,t],
 \eeqlb
which is given by
 \beqlb\label{eq2.13}
\pi_{r,t}(\lambda)= \lambda\prod_{s\in(r,t]} (1+\Delta\alpha(s)) \exp\{\alpha_c(t)-\alpha_c(r)\}.
 \eeqlb

 \eitemize
\ecorollary

We next discuss briefly the structures of the transition semigroup $\{Q_{r,t}: t\ge r\in I\}$ defined by \eqref{eq1.1}. A family of mappings $\{v_{r,t}: t\ge r\in I\}$ on $(0,\infty)$ is called a \textit{cumulant semigroup} if the following conditions are satisfied:
 \bitemize

\itm[(2.A)] (Semigroup property) for $\lambda> 0$ and $t\ge s\ge r\in I$,
 \beqlb\label{eq2.14}
v_{r,t}(\lambda)= v_{r,s}\circ v_{s,t}(\lambda)= v_{r,s}(v_{s,t}(\lambda));
 \eeqlb

\itm[(2.B)] (L\'{e}vy--Kthintchine representation) for $\lambda> 0$ and $t\ge r\in I$,
 \beqlb\label{eq2.15}
v_{r,t}(\lambda)= a_{r,t} + h_{r,t}\lambda + \int_0^\infty(1-\e^{-\lambda y}) l_{r,t}(\d y),
 \eeqlb
where $a_{r,t}\ge 0$, $h_{r,t}\ge 0$ and $(1\wedge y) l_{r,t}(\d y)$ is a finite measure on $(0,\infty)$.

 \eitemize
Given a cumulant semigroup $\{v_{r,t}: t\ge r\in I\}$, we can define the transition semigroup $\{Q_{r,t}: t\ge r\in I\}$ on $[0,\infty]$ using \eqref{eq1.1}. Clearly, the CBVE-process with this transition semigroup has both $0$ and $\infty$ as traps. By \eqref{eq1.1} we have
 \beqlb\label{eq2.16}
Q_{r,t}(x,[0,\infty))= \e^{-xv_{r,t}(0)}, \qquad x\in [0,\infty),
 \eeqlb
and
 \beqlb\label{eq2.17}
Q_{r,t}(x,\{0\})= \e^{-xv_{r,t}(\infty)}, \qquad x\in (0,\infty),
 \eeqlb
where $v_{r,t}(0):= \lim_{\lambda\downarrow 0} v_{r,t}(\lambda)= a_{r,t}\in [0,\infty)$ and $v_{r,t}(\infty):= \lim_{\lambda\uparrow \infty} v_{r,t}(\lambda)\in (0,\infty]$. We say a cumulant semigroup $\{v_{r,t}: t\ge r\in I\}$ is \textit{conservative} if $v_{r,t}(0)= a_{r,t}= 0$ for all $t\ge r\in I$. In this case, we can restrict $\{Q_{r,t}: t\ge r\in I\}$ to a conservative transition semigroup on $[0,\infty)$ and rewrite \eqref{eq1.1} into
 \beqlb\label{eq2.18}
\int_{[0,\infty)}\e^{-\lambda y}Q_{r,t}(x,\d y)= \e^{-xv_{r,t}(\lambda)}, \qquad x\ge 0,\lambda\ge 0.
 \eeqlb

To conclude this section, we prove some useful upper and lower bounds for the solutions to the integral evolution equation \eqref{eq1.3}. Let $(b_1,c,m)$ be admissible parameters. For $\lambda> 0$ and $t\ge r\ge 0$ let
 \beqlb\label{eq2.19}
U_{r,t}(\lambda)= [\lambda + m((0,t]\times (1,\infty))] \exp\{\|b_1\|(t)-\|b_1\|(r)\}.
 \eeqlb
By the admissibility of the parameters we have $m(\{s\}\times (0,1])= 0$ and $m(\{s\}\times (1,\infty))> 0$ for $s\in J$. For $t\ge 0$ choose a sufficiently large constant $\eta_t>1$ so that $m(\{s\}\times (1,\eta_t])> 0$ when $s\in (0,t]\cap J$. Let
 \beqnn
F_t(\lambda)= U_{0,t}(\lambda)^{-1}(1-\e^{-U_{0,t}(\lambda)}),
 \quad
H_t(\lambda)= [\eta_tU_{0,t}(\lambda)]^{-1}(1-\e^{-\eta_tU_{0,t}(\lambda)}).
 \eeqnn
Let $\alpha(r)= \alpha(r,t,\lambda)$ be the c\`{a}dl\`{a}g function on $[0,t]$ defined by
 \beqlb\label{eq2.20}
\alpha(r)\ar=\ar - \frac{1}{2}U_{0,t}(\lambda) \int_0^r\int_0^{\varepsilon_t(\lambda)} z^2 m(\d s,\d z) + H_t(\lambda) \int_0^r\int_1^{\eta_t} z m(\d s,\d z) \cr
 \ar\ar
-\, b_1(r) - U_{0,t}(\lambda)c(r) - [1-F_t(\lambda)]\int_0^r\int_{\varepsilon_t(\lambda)}^1 z m(\d s,\d z),
 \eeqlb
where $\varepsilon_t(\lambda)= 1\land [U_{0,t}(\lambda)^{-1}F_t(\lambda)]$. Let
 \beqlb\label{eq2.21}
l_{r,t}(\lambda)= \lambda\prod_{s\in(r,t]} [1+(0\land\Delta\alpha(s))] \exp\{\|\alpha\|(r)-\|\alpha\|(t)\}.
 \eeqlb

\bproposition\label{th2.4} Suppose that $r\mapsto v_{r,t}(\lambda)$ is a bounded positive solution to \eqref{eq1.3} with $\lambda> 0$. Then we have
 \beqlb\label{eq2.22}
l_{r,t}(\lambda)\le v_{r,t}(\lambda)\le U_{r,t}(\lambda), \qquad r\in [0,t].
 \eeqlb
\eproposition

\proof The upper bound in \eqref{eq2.22} follows by Gronwall's inequality since \eqref{eq1.3} implies
 \beqnn
v_{r,t}(\lambda)\le \lambda + \int_0^t\int_1^\infty m(\d s,\d z) + \int_r^t v_{s,t}(\lambda)\|b_1\|(\d s).
 \eeqnn
Let $r\mapsto \pi_{r,t}(\lambda)$ be the solution to \eqref{eq2.12} with $\alpha$ given by \eqref{eq2.20}. Then we have
 \beqlb\label{eq2.23}
v_{r,t}(\lambda) - \pi_{r,t}(\lambda)
 =
G_{r,t}(\lambda) + \int_r^t [v_{s,t}(\lambda) - \pi_{s,t}(\lambda)] \alpha(\d s),
 \eeqlb
where
 \beqnn
G_{r,t}(\lambda)\ar=\ar \int_r^t v_{s,t}(\lambda)[U_{0,t}(\lambda) - v_{s,t}(\lambda)] c(\d s) + \int_r^t\int_{\eta_t}^\infty \big(1 - \e^{-v_{s,t} (\lambda)z}\big)m(\d s,\d z) \cr
 \ar\ar\quad
+ \int_r^t\int_0^{\varepsilon_t(\lambda)}\Big[\frac{1}{2}U_{0,t}(\lambda)v_{s,t}(\lambda)z^2 - K(v_{s,t}(\lambda),z)\Big] m(\d s,\d z) \cr
 \ar\ar\quad
+ \int_r^t\int_{\varepsilon_t(\lambda)}^1 \big[1-\e^{-v_{s,t}(\lambda)z} - F_t(\lambda)v_{s,t}(\lambda)z\big] m(\d s,\d z) \cr
 \ar\ar\quad
+ \int_r^t\int_1^{\eta_t} \big[1-\e^{-v_{s,t}(\lambda)z} - H_t(\lambda)v_{s,t}(\lambda)z\big] m(\d s,\d z),
 \eeqnn
where $K(\lambda,z)= \e^{-\lambda z} - 1 + \lambda z$. In view of \eqref{eq1.4}, for any $s\in(0,t]$ we have
 \beqnn
\int_0^1 z m_d(\{s\},\d z)
 \le
1-\Delta b_1(s).
 \eeqnn
It follows that
 \beqnn
\Delta\alpha(s)\ar=\ar - \frac{1}{2}U_{0,t}(\lambda) \int_0^{\varepsilon_t(\lambda)} z^2 m_d(\{s\},\d z) - [1-F_t(\lambda)] \int_{\varepsilon_t(\lambda)}^1 z m_d(\{s\},\d z) \cr
 \ar\ar
-\, \Delta b_1(s) + H_t(\lambda) \int_1^{\eta_t} z m(\{s\},\d z) \cr
 \ar\ge\ar
- \frac{1}{2}U_{0,t}(\lambda){\varepsilon_t(\lambda)} \int_0^1 z m_d(\{s\},\d z) - [1-F_t(\lambda)] \int_0^1 z m_d(\{s\},\d z) \cr
 \ar\ar
-\,\Delta b_1(s) + H_t(\lambda) \int_1^{\eta_t} z m(\{s\},\d z) \cr
 \ar\ge\ar
- \frac{1}{2}U_{0,t}(\lambda){\varepsilon_t(\lambda)} [1-\Delta b_1(s)] - [1-F_t(\lambda)] [1-\Delta b_1(s)] \cr
 \ar\ar
-\,\Delta b_1(s) + H_t(\lambda) \int_1^{\eta_t} z m(\{s\},\d z).
 \eeqnn
By the admissibility of the parameters we have $1-\Delta b_1(s)> 0$ when $s\in (0,t]\setminus J$ and $m(\{s\}\times (1,\eta_t])> 0$ when $s\in (0,t]\cap J$, so $\Delta\alpha(s)>-1$ for each $s\in (0,t]$. Then Proposition~\ref{th2.2} applies to \eqref{eq2.23}. Since $r\mapsto G_{r,t}(\lambda)$ is a decreasing function, from \eqref{eq2.9} we see $v_{r,t}(\lambda) - \pi_{r,t}(\lambda)\ge 0$. By comparing \eqref{eq2.13} and \eqref{eq2.21} we have the lower bound in \eqref{eq2.22}. \qed


\section{Conservative cumulant semigroups}

\setcounter{equation}{0}

In this section, we take $I=[0,\infty)$. Let $\alpha$ be a c\`{a}dl\`{a}g function on $[0,\infty)$ having locally bounded variations and satisfying $\Delta\alpha(t)>-1$ for $t> 0$. Let $\mu(\d s,\d z)$ be a $\sigma$-finite measure on $(0,\infty)^2$ satisfying
 \beqlb\label{eq3.1}
\int_0^t\int_0^\infty z\mu(\d s,\d z)<\infty, \qquad t\ge 0.
 \eeqlb
Given $t\ge 0$ and $\lambda\ge 0$, we first consider the backward integral evolution equation:
 \beqlb\label{eq3.2}
u_{r,t}(\lambda)= \lambda + \int_r^t u_{s,t}(\lambda)\alpha(\d s) + \int_r^t\int_0^\infty\big(1 - \e^{-u_{s,t}(\lambda) z}\big)\mu(\d s,\d z), \qquad r\in [0,t].
 \eeqlb
This is clearly spacial case of \eqref{eq1.3}.

\bproposition\label{th3.1} For $t\ge 0$ and $\lambda\ge 0$, there is a unique bounded positive solution $r\mapsto u_{r,t}(\lambda)$ on $[0,t]$ to \eqref{eq3.2} and $(u_{r,t})_{t\ge r}$ is a conservative cumulant semigroup. Moreover, for $\lambda\ge 0$ we have
 \beqlb\label{eq3.3}
u_{r,t}(\lambda)\le \lambda\e^{\|\rho\|(r,t]}
\le
\lambda\e^{\|\rho\|(t)},
 \eeqlb
where
 \beqnn
\rho(t)= \alpha (t) + \int_0^t\int_0^\infty z \mu(\d s,\d z).
 \eeqnn
\eproposition

\proof \textit{Step~1}. Let $r\mapsto u_{r,t}(\lambda)$ be a bounded positive solutions to \eqref{eq3.2}. From the equation it is easy to see that
 \beqnn
u_{r,t}(\lambda)
 \le
\lambda + \int_r^t u_{s,t}(\lambda) \|\rho\|(\d s).
 \eeqnn
Then \eqref{eq3.3} follows by Gronwall's inequality. Suppose that $r\mapsto w_{r,t}(\lambda)$ is also a bounded positive solution to \eqref{eq3.2}. Then we have
 \beqnn
|u_{r,t}(\lambda) - w_{r,t}(\lambda)|
 \le
\int_r^t|u_{s,t}(\lambda) - w_{s,t}(\lambda)|\|\rho\|(\d s).
 \eeqnn
By Gronwall's inequality we see $|u_{r,t}(\lambda) - w_{r,t}(\lambda)|= 0$ for every $r\in [0,t]$.

\textit{Step~2}. Consider the case where $\alpha$ vanishes. Let $t\ge 0$ and $\lambda\ge 0$ be fixed. For $r\in [0,t]$ set $v^{(0)}_{r,t}(\lambda)= 0$ and define $v^{(k)}_{r,t}(\lambda)$ inductively by
 \beqlb\label{eq3.4}
v^{(k+1)}_{r,t}(\lambda)\ar=\ar \lambda + \int_r^t\int_0^\infty \big(1 - \e^{-v^{(k)}_{s,t}(\lambda) z}\big)\mu(\d s,\d z).
 \eeqlb
By Proposition~4.2 of Silverstein (1968) one can use \eqref{eq3.4} to see inductively that each $v^{(k)}_{r,t}(\lambda)$ has the L\'{e}vy--Kthintchine representation \eqref{eq2.15}. Moreover, we have $v^{(k)}_{r,t}(0)= 0$ and
 \beqnn
0\le v^{(k)}_{r,t}(\lambda)\le v^{(k+1)}_{r,t}(\lambda)\le \pi_{r,t}(\lambda),
 \eeqnn
where $r\mapsto \pi_{r,t}(\lambda)$ is the solution to \eqref{eq2.12} with
 \beqnn
\alpha(t)= \int_0^t\int_0^\infty z \mu(\d s,\d z).
 \eeqnn
Then the limit $v_{r,t}(\lambda)= \uparrow\lim_{k\to \infty} v^{(k)}_{r,t}(\lambda)$ exists and the convergence is uniform in $(r,\lambda)\in [0,t]\times [0,B]$ for every $t\ge 0$ and $B\ge 0$. In fact, setting $u_k (r,t,\lambda)= \sup_{r\le s\le t}|v^{(k)}_{s,t}(\lambda) - v^{(k-1)}_{s,t}(\lambda)|$, we have
 \beqnn
u_k (r,t,\lambda)
 \ar\le\ar
\int_r^t u_{k-1}(t_1,t,\lambda) \alpha(\d t_1) \cr
 \ar\le\ar
\int_r^t\alpha(\d t_1)\int_{t_1}^t u_{k-2} (t_2,t,\lambda) \alpha(\d t_2)\le \cdots \cr
 \ar\le\ar
\int_r^t \alpha(\d t_1)\int_{t_1}^t\cdots\int_{t_{k-2}}^t u_1 (t_{k-1},t,\lambda) \alpha(\d t_{k-1}) \cr
 \ar\le\ar
B\int_r^t\alpha(\d t_1)\int_{t_1}^t\cdots\int_{t_{k-2}}^t \alpha(\d t_{k-1}) \cr
 \ar\le\ar
B\frac{\alpha(r,t]^{k-1}}{(k-1)!}
 \le
B\frac{\alpha(0,t]^{k-1}}{(k-1)!},
 \eeqnn
and hence $\sum_{k = 1}^\infty u_k(r,t,\lambda)\le B\e^{\alpha(0,t]}< \infty$. By Corollary~1.33 in Li (2011) we infer that $v_{r,t}(\lambda)$ has representation \eqref{eq2.15} with $v_{r,t}(0)= 0$. By \eqref{eq3.4} and monotone convergence we see $r\mapsto v_{r,t}(\lambda)$ is a solution to \eqref{eq3.2} with $\alpha\equiv 0$. The semigroup property \eqref{eq2.14} follows from the uniqueness of the solution. Then $(u_{r,t})_{t\ge r}$ is a conservative cumulant semigroup.

\textit{Step~3}. Let $\zeta$ be the c\`{a}dl\`{a}g function on $[0,\infty)$ such that $\zeta_c(t)= \alpha_c(t)$ and $\Delta\zeta(t)= \log [1+\Delta\alpha(t)]$ for every $t\ge 0$. By the second step, there is a unique bounded positive solution $r\mapsto u_{r,t}(\lambda)$ on $[0,t]$ to
 \beqlb\label{eq3.5}
u_{r,t}(\lambda)= \lambda + \int_r^t\int_0^\infty (1-\e^{-u_{s,t}(\lambda)z}) \e^{\zeta(s-)}\mu(\d s,\e^{\zeta(s)}\d z).
 \eeqlb
Moreover, the family $(u_{r,t})_{t\ge r}$ is a conservative cumulant semigroup. Then we can define another conservative cumulant semigroup $(v_{r,t})_{t\ge r}$ by $v_{r,t}(\lambda)= \e^{-\zeta(r)} u_{r,t}(\e^{\zeta(t)}\lambda)$. Since $u_{t,t}(\lambda)= \lambda$, by integration by parts we have
 \beqnn
\lambda\ar=\ar \e^{-\zeta(r)}u_{r,t}(\e^{\zeta(t)}\lambda) + \int_r^t u_{s,t}(\e^{\zeta(t)}\lambda) \d\e^{-\zeta(s)} + \int_r^t \e^{-\zeta(s-)} \d u_{s,t}(\e^{\zeta(t)}\lambda) \cr
 \ar=\ar
v_{r,t}(\lambda) - \int_r^t v_{s,t}(\lambda) \alpha(\d s) - \int_r^t\int_0^\infty (1-\e^{-v_{s,t}(\lambda)z}) \mu(\d s,\d z).
 \eeqnn
Then $r\mapsto v_{r,t}(\lambda)$ is a solution to \eqref{eq3.2} on $[0,t]$. \qed

We next consider a more interesting spacial case of \eqref{eq1.3}. Let $(b_1,c,m)$ be admissible parameters given as in the introduction. Instead of \eqref{eq1.2}, we here assume the stronger integrability condition:
 \beqlb\label{eq3.6}
m(t):= \int_0^t\int_0^\infty(z\wedge z^2)m(\d s,\d z)< \infty, \qquad t\ge 0.
 \eeqlb
Then we can rewrite \eqref{eq1.3} equivalently into:
 \beqlb\label{eq3.7}
v_{r,t}(\lambda)= \lambda - \int_r^t v_{s,t}(\lambda)b(\d s) - \int_r^t v_{s,t}(\lambda)^2c(\d s) - \int_r^t\int_0^\infty K(v_{s,t}(\lambda),z) m(\d s,\d z),
 \eeqlb
where $K(\lambda,z)= \e^{-\lambda z} - 1 + \lambda z$ and
 \beqlb\label{eq3.8}
b(t)= b_1(t) - \int_0^t\int_1^\infty z m(\d s,\d z).
 \eeqlb

Let $B[0,\infty)^+$ be the set of locally bounded positive Borel functions on $[0,\infty)$ and $M[0,\infty)$ the set of Radon measures on $[0,\infty)$. By a \textit{branching mechanism} with parameters $(b,c,m)$ we mean the functional $\phi$ on $B[0,\infty)^+\times M[0,\infty)$ defined by
 \beqnn
\phi(f,B)
 \ar=\ar
\int_{B}f(s)b(\d s) + \int_{B}f(s)^2c(\d s) + \int_B\int_0^\infty K(f(s),z) m(\d s,\d z),
 \eeqnn
where $f\in B[0,\infty)^+$ and $B\in\mcr{B}[0,\infty)$. Using this notation, we can rewrite \eqref{eq3.7} equivalently into
 \beqlb\label{eq3.9}
v_{r,t}(\lambda)= \lambda - \phi(v_{\cdot,t}(\lambda),(r,t]), \qquad r\in [0,t].
 \eeqlb
For any integer $n\ge 1$ we define the branching mechanism $\phi_n$ on by
 \beqlb\label{eq3.10}
\phi_n(f,B)
 \ar=\ar
\int_{B}f(s)b(\d s) + 2n^2\int_{B}\big(\e^{-f(s)/n} - 1 + f(s)/n\big)c(\d s) \cr
 \ar\ar
+ \int_{B}\int_0^1 \big(\e^{-f(s)z} - 1 + f(s)z\big)(1\wedge(nz))m(\d s,\d z) \cr
 \ar\ar
+ \int_B\int_1^\infty \big(\e^{-f(s)z} - 1 + f(s)z\big)m(\d s,\d z) \cr
 \ar=\ar
-\int_{B}f(s)\alpha_n(\d s) - \int_{B}\int_0^\infty\big(1 - \e^{-f(s)z}\big)\mu_n(\d s,\d z),
 \eeqlb
where
 \beqnn
\alpha_n(\d s)\ar=\ar -b(\d s) - 2nc(\d s) - \int_0^1 z(1\wedge nz) m(\d s,\d z) \cr
 \ar\ar\qqquad
- \int_1^\infty z m(\d s,\d z) \cr
 \ar=\ar
-b_1(\d s) - 2nc(\d s) - \int_0^1 z(1\wedge nz) m(\d s,\d z)
 \eeqnn
and
 \beqnn
\mu_n(\d s,\d z)\ar=\ar 2n^2c(\d s)\delta_{1/n}(\d z) + 1_{\{z\le 1\}} (1\wedge nz) m(\d s,\d z) \ccr
 \ar\ar\qqquad
+\, 1_{\{z>1\}} m(\d s,\d z).
 \eeqnn
Then $\Delta\alpha_n(t)>-1$ for every $t> 0$ since $(b_1,c,m)$ are admissible parameters.

\blemma\label{th3.2} The branching mechanisms $\phi$ and $\phi_n$ have the following properties:
 \bitemize

\itm[{\rm(i)}] For $t\ge r\ge 0$ and $f\in B[0,\infty)^+$, we have $\phi(f,(r,t])= \uparrow \lim_{n\uparrow\infty}\phi_n(f,(r,t])$;

\itm[{\rm(ii)}] For $t\ge s\ge r\ge 0$ and $f\le g\in B[0,\infty)^+$,
 \beqnn
\phi(f,(s,t]) - \phi_n(f,(s,t])\le \phi(g,(r,t]) - \phi_n(g,(r,t]);
 \eeqnn

\itm[{\rm(iii)}] For $t\ge s\ge r\ge 0$ and $f, g\in B[0,\infty)^+$,
 \beqnn
\big|\phi(f,(s,t])-\phi(g,(s,t])\big|
 \le
\big[C_1(t)+1\big]\int_r^t|f(s) - g(s)|C_2(\d s),
 \eeqnn
where $C_1(t)= \sup_{s\in[0,t]}[f(s) + g(s)]$ and
 \beqlb\label{eq3.11}
C_2(\d s)= \|b\|(\d s) + c(\d s) + \int_0^\infty(z\wedge z^2)m(\d s,\d z).
 \eeqlb
 \eitemize
\elemma

\proof By \eqref{eq3.10} we obtain immediately (i) and (ii). For any $t\ge s\ge r\ge 0$ and $f, g\in B[0,\infty)^+$, we have
 \beqnn
|\phi(f,(s,t]) - \phi(g,(s,t])|
 \ar\le\ar
\int_r^t|f(u) - g(u)|\|b\|(\d u) + C_1(t)\int_r^t|f(u) - g(u)|c(\d u) \cr
 \ar\ar\qquad
+\, C_1(t)\int_r^t\int_0^1|f(u) - g(u)|z^2m(\d u,\d z) \cr
 \ar\ar\qquad
+\, \int_r^t\int_1^\infty|f(u) - g(u)|z m(\d u,\d z).
 \eeqnn
Then $\mrm{(iii)}$ follows.\qed

\btheorem\label{th3.3} For every $t\ge 0$ and $\lambda\ge0$ there is a unique bounded positive solution $r\mapsto v_{r,t}(\lambda)$ to the integral evolution equation \eqref{eq3.7} or \eqref{eq3.9} and $(v_{r,t})_{t\ge r}$ is a conservative cumulant semigroup.
\etheorem

\proof Let $\phi_n$ be defined by \eqref{eq3.10}. It is easy to see that $\alpha_n$ and $\mu_n$ satisfy the conditions of Proposition~\ref{th3.1}. In particular, for any $s>0$ we have
 \beqnn
\Delta\alpha_n(s)\ge -\Delta b(s) - (1-\e^{-n})\int_0^\infty z m(\{s\},\d z)> -1.
 \eeqnn
Then a conservative cumulant semigroup $(v^{(n)}_{r,t})_{t\ge r}$ is defined by the evolution integral equation:
 \beqlb\label{eq3.12}
v^{(n)}_{r,t}(\lambda)= \lambda - \phi_n(v^{(n)}_{\cdot,t}(\lambda),(r,t]), \qquad \lambda\ge 0, r\in [0,t].
 \eeqlb
By \eqref{eq3.3} we have $v^{(n)}_{r,t}(\lambda)\le A\e^{\|b\|(t)}$ for $r\in [0,t]$ and $\lambda\in [0,A]$. For $n\ge k\ge 1$ let
 \beqnn
D_{k,n}(r,t,\lambda)= \sup_{r\le s\le t}\big|v^{(n)}_{s,t}(\lambda) - v^{(k)}_{s,t}(\lambda)\big|.
 \eeqnn
By Lemma~\ref{th3.2}, we have
 \beqnn
D_{k,n}(r,t,\lambda)
 \ar\le\ar
2|\phi(A\e^{\|b\|(t)},(0,t]) - \phi_k(A\e^{\|b\|(t)},(0,t])| \ccr
 \ar\ar
+\, \big[C_1(t)+1\big]\int_r^tD_{k,n}(s,t,\lambda) C_2(\d s),
 \eeqnn
where $C_1(t)= 2A\e^{\|b\|(t)}$ and $C_2(\d s)$ is given by \eqref{eq3.11}. By Gronwall's inequality,
 \beqnn
D_{k,n} (r,t,\lambda)
 \le
2|\phi(A\e^{\|b\|(t)},(0,t])- \phi_k(A\e^{\|b\|(t)},(0,t])|\e^{[C_1(t)+1]C_2(t)}.
 \eeqnn
By Lemma~\ref{th3.2} it is easy to see the limit $v_{r,t}(\lambda):= \lim_{k\to \infty} v^{(k)}_{r,t}(\lambda)$ exists and convergence is uniform in $(r,\lambda)\in [0,t] \times [0,A]$ for every $A\ge 0$. By Corollary~1.33 in Li (2011) we have the L\'{e}vy--Kthintchine representation \eqref{eq2.15} with $v_{r,t}(0)= 0$. From \eqref{eq3.12} we get \eqref{eq3.7} by Lemma~\ref{th3.2} and dominated convergence. The uniqueness of bounded positive solution to \eqref{eq3.7} follows by Lemma~\ref{th3.2} and Gronwall's inequality. The semigroup property \eqref{eq2.14} follows from the uniqueness of the solution. Then $(v_{r,t})_{t\ge r}$ is a conservative cumulant semigroup. \qed

\bproposition\label{th3.4} Let $r\mapsto v_{r,t}(\lambda)$ be the unique bounded positive solution to \eqref{eq3.7} on $[0,t]$. Then we have $v_{r,t}(\lambda)\le \pi_{r,t}(\lambda)$, where
 \beqlb\label{eq3.13}
r\mapsto \pi_{r,t}(\lambda):= \lambda\prod_{r<s\le t}(1 - \Delta b(s))\exp\{b_c(r)- b_c(t)\}
 \eeqlb
is the solution to \eqref{eq2.12} with $\alpha(t)= -b(t)$.
\eproposition

\proof Fix $t\ge 0$ and $\lambda\ge 0$ and let $H(r)= \pi_{r,t}(\lambda)- v_{r,t}(\lambda)$. From \eqref{eq2.12} and \eqref{eq3.7} we see that $r\mapsto H(r)$ satisfies \eqref{eq2.8} with $\alpha(t)= -b(t)$ and
 \beqnn
G(r)= \int_r^tv_{s,t}(\lambda)^2 c(\d s) + \int_r^t\int_0^\infty K(v_{s,t}(\lambda),z) m(\d s,\d z).
 \eeqnn
By Proposition~\ref{th2.2} we have the representation \eqref{eq2.9} for $H(r)$, which implies $H(r)\ge 0$ since $r\mapsto G(r)$ is a decreasing function on $[0,t]$. \qed

\bproposition\label{th3.5} For $t\ge r\ge 0$ let $\pi_{r,t}(1)$ be defined by \eqref{eq3.13} with $\lambda=1$. Then we have
 \beqlb\label{eq3.14}
h_{r,t} + \int_0^\infty yl_{r,t}(\d y)= \pi_{r,t}(1)
 \eeqlb
and
 \beqlb\label{eq3.15}
\int_{[0,\infty)}y Q_{r,t}(x,\d y)= x\pi_{r,t}(1), \qquad x\ge 0.
 \eeqlb
\eproposition

\proof By Proposition~\ref{th3.4} we have $\lambda^{-1} v_{r,t}(\lambda)\le \pi_{r,t}(1)$. Then we can differentiate both sides of \eqref{eq3.7} and use bounded convergence to see that $r\mapsto \frac{\partial}{\partial \lambda}v_{r,t}(0+)$ is a solution to \eqref{eq2.12} with $\alpha(t)= -b(t)$ and $\lambda=1$. It follows that $\frac{\partial}{\partial \lambda}v_{r,t}(0+)\equiv \pi_{r,t}(1)$. By differentiating both sides of \eqref{eq2.15} we obtain \eqref{eq3.14}. Similarly we get \eqref{eq3.15} from \eqref{eq2.18}. \qed

The transformation of the cumulant semigroup used in the proof of Proposition~\ref{th3.1} is an inhomogeneous nonlinear variation of the classical $h$-transformation and has been used in the study of CB-processes; see, e.g., Bansaye et al.\ (2013), He et al.\ (2018) and Li (2011, Section~6.1). A generalized form of the transformation is given below, which will be useful in the next section.

\btheorem\label{th3.6} Let $(v_{r,t})_{t\ge r}$ be the conservative cumulant semigroup defined by \eqref{eq3.7} or \eqref{eq3.9}. Let $t\mapsto \zeta(t)$ be a locally bounded function on $[0,\infty)$. Then another conservative cumulant semigroup $(u_{r,t})_{t\ge r}$ is defined by:
 \beqlb\label{eq3.16}
u_{r,t}(\lambda) = \e^{\zeta(r)} v_{r,t}(\e^{-\zeta(t)}\lambda),\qquad \lambda\ge 0.
 \eeqlb
Moreover, if $\zeta$ is a c\`{a}dl\`{a}g function on $[0,\infty)$ with locally bounded variations, then $[0,t]\ni r\mapsto u_{r,t}(\lambda)$ is the unique bounded positive solution to
 \beqlb\label{eq3.17}
u_{r,t}(\lambda)\ar=\ar \lambda - \int_r^t u_{s,t}(\lambda) \d\beta(s) - \int_r^t u_{s,t}(\lambda) \e^{-\Delta\zeta(s)}b(\d s) - \int_r^t u_{s,t}(\lambda)^2 \e^{-\zeta(s)} c(\d s) \cr
 \ar\ar\qqquad
- \int_r^t\int_0^\infty K(u_{s,t}(\lambda),z) \e^{\zeta(s-)}m(\d s,\e^{\zeta(s)}\d z),
 \eeqlb
where
 \beqnn
\beta(t)= \zeta_c(t) + \sum_{s\in (0,t]} \big(1-\e^{-\Delta\zeta(s)}\big).
 \eeqnn
\etheorem

\proof The arguments are generalizations of those in the last step of the proof of Proposition~\ref{th3.1}. Clearly, the family $(u_{r,t})_{t\ge r}$ defined by \eqref{eq3.16} is a conservative cumulant semigroup. If $\zeta$ is a c\`{a}dl\`{a}g function with locally bounded variations, we can use integration by parts to get
 \beqnn
\lambda\ar=\ar \e^{\zeta(r)}v_{r,t}(\e^{-\zeta(t)}\lambda) + \int_r^t v_{s,t}(\e^{-\zeta(t)}\lambda) \d\e^{\zeta(s)} + \int_r^t \e^{\zeta(s-)} \d v_{s,t}(\e^{-\zeta(t)}\lambda) \cr
 \ar=\ar
u_{r,t}(\lambda) + \int_r^t v_{s,t}(\e^{-\zeta(t)}\lambda)\e^{\zeta(s)} \zeta_c(\d s) + \sum_{s\in(0,t]} v_{s,t}(\e^{-\zeta(t)}\lambda)(\e^{\zeta(s)}-\e^{\zeta(s-)}) \cr
 \ar\ar
+ \int_r^t \e^{\zeta(s-)}v_{s,t}(\e^{-\zeta(t)}\lambda) b(\d s) + \int_r^t \e^{\zeta(s-)}v_{s,t}(\e^{-\zeta(t)}\lambda)^2 c(\d s) \cr
 \ar\ar
+ \int_r^t\int_0^\infty \e^{\zeta(s-)}K(v_{s,t}(\e^{-\zeta(t)}\lambda),z) m(\d s,\d z) \cr
 \ar=\ar
u_{r,t}(\lambda) + \int_r^t v_{s,t}(\e^{-\zeta(t)}\lambda)\e^{\zeta(s)} \beta(\d s) + \int_r^t \e^{\zeta(s-)}v_{s,t}(\e^{-\zeta(t)}\lambda) b(\d s) \cr
 \ar\ar
+ \int_r^t \e^{\zeta(s-)}v_{s,t}(\e^{-\zeta(t)}\lambda)^2 c(\d s) + \int_r^t\int_0^\infty \e^{\zeta(s-)}K(v_{s,t}(\e^{-\zeta(t)}\lambda),z) m(\d s,\d z) \cr
 \ar=\ar
u_{r,t}(\lambda) + \int_r^t u_{s,t}(\lambda) \beta(\d s) + \int_r^t \e^{-\Delta\zeta(s)} u_{s,t}(\lambda) b(\d s) + \int_r^t \e^{-\zeta(s)-\Delta\zeta(s)}u_{s,t}(\lambda)^2 c(\d s) \cr
 \ar\ar
+ \int_r^t\int_0^\infty \e^{\zeta(s-)} K(u_{s,t}(\lambda),\e^{-\zeta(s)}z) m(\d s,\d z) \cr
 \ar=\ar
u_{r,t}(\lambda) + \int_r^t u_{s,t}(\lambda) \beta(\d s) + \int_r^t \e^{-\Delta\zeta(s)} u_{s,t}(\lambda) b(\d s) + \int_r^t \e^{-\zeta(s)}u_{s,t}(\lambda)^2 c(\d s) \cr
 \ar\ar
+ \int_r^t\int_0^\infty \e^{\zeta(s-)} K(u_{s,t}(\lambda),z) m(\d s,\e^{\zeta(s)}\d z),
 \eeqnn
where we have used the continuity of $s\mapsto c(s)$ for the last equality. Then $r\mapsto u_{r,t}(\lambda)$ solves \eqref{eq3.17}. The uniqueness of the solution holds by Theorem~\ref{th3.3}. \qed


\section{Stochastic equations for CBVE-processes}

\setcounter{equation}{0}

Let $(b,c,m)$ be given as in the last section. Recall that $m$ also denotes the increasing function defined by \eqref{eq3.6}. Then $J_m:= \{s>0: \Delta m(s)> 0\}$ is at most a countable set. Let $m_d(\d s,\d z)= 1_{J_m}(s)m(\d s,\d z)$ and $m_c(\d s,\d z)= m(\d s,\d z) - m_d(\d s,\d z)$. Suppose that $(\Omega, \mcr{F}, \mcr{F}_t, \mbf{P})$ is a filtered probability space satisfying the usual hypotheses. Let $W(\d s,\d u)$ and $M(\d s,\d z,\d u)$ be $(\mcr{F}_t)$-noises given as in the introduction. One can see that $M_c(\d s,\d z,\d u):= 1_{J_m^c}(s) M(\d s,\d z,\d u)$ and $M_d(\d s,\d z,\d u):= 1_{J_m}(s) M(\d s,\d z,\d u)$ are $(\mcr{F}_t)$-Poisson random measures with intensities $m_c(\d s,\d z)\d u$ and $m_d(\d s,\d z)\d u$, respectively. Those random measures are independent of each other as they have disjoint supports. We can rewrite \eqref{eq1.5} equivalently into:
 \beqlb\label{eq4.1}
X(t)\ar=\ar X(0) + \int_0^t\int_0^{X(s-)}W(\d s,\d u) + \int_0^t\int_0^\infty \int_0^{X(s-)} z\tilde{M}_c(\d s,\d z,\d u) \cr
 \ar\ar\qquad
- \int_0^tX(s-)b(\d s) + \int_0^t\int_0^\infty\int_0^{X(s-)}z\tilde{M}_d(\d s,\d z,\d u),
 \eeqlb
where $b$ is defined by \eqref{eq3.8}.

\bproposition\label{th4.1} Let $\{X(t): t\ge 0\}$ be a solution to \eqref{eq4.1} and let $\tau_k = \inf\{t\ge0: X(t)\ge k\}$ for $k\ge 1$. Then $\tau_k\to \infty$ almost surely as $k\to\infty$. Moreover, for $t\ge 0$ and $k\ge 1$ we have
 \beqlb\label{eq4.2}
k\mbf{P}\{\tau_k\le t\}\le \mbf{P}[X(0)] \e^{\|b\|(t)}
 \quad\mbox{and}\quad
\mbf{P}[X(t)]\le \mbf{P}[X(0)] \e^{\|b\|(t)}.
 \eeqlb
\eproposition

 %
 We omit the proof of the above proposition, which is based on an application of Gronwall's inequality.
The comparison property of the solutions to \eqref{eq4.1} plays an important role in the analysis of the stochastic equation. In the proof, some special care has to be taken for the negative jumps brought about by the compensator of the Poisson random measure. For simplicity we only give a treatment of the property under a stronger integrability condition.

\bproposition\label{th4.2} The pathwise uniqueness of solution holds for \eqref{eq4.1} under the additional integrability condition
 \beqlb\label{eq4.4}
\int_0^t\int_0^\infty z^2 m(\d s,\d z)< \infty, \qquad t\ge 0.
 \eeqlb
Moreover, under the above condition, if $\{X_1(t): t\ge0\}$ and $\{X_2(t): t\ge0\}$ are two solutions to \eqref{eq4.1} satisfying $\mbf{P}\{X_1(0)\le X_2(0)\}= 1$, then we have $\mbf{P}\{X_1(t)\le X_2(t)$ for every $t\ge 0\}= 1$.
\eproposition

\proof It suffices to prove the second assertion. For each integer $n\ge 0$ define $a_n= \exp\{-n(n+1)/2\}$. Then $\int_{a_n}^{a_{n-1}}z^{-1} \d z = n$ and $a_n\to 0$ decreasingly as $n\to \infty$. Let $x\mapsto g_n(x)$ be a positive continuous function supported by $(a_n,a_{n-1})$ so that $\int_{a_n}^{a_{n-1}}g_n(x)\d x=1$ and $g_n(x)\le 2(nx)^{-1}$ for every $x>0$. For $n\ge 0$ and $z\in \mbb{R}$ let
 \beqlb\label{eq4.5}
f_n(z)=\int_0^{z\vee 0}\d y\int_0^yg_n(x)\d x.
 \eeqlb
From \eqref{eq4.1} we have
 \beqlb\label{eq4.6}
X(t)\ar=\ar X(0) + \int_0^t\int_0^{X(s-)}W(\d s,\d u) + \int_0^t\int_0^\infty \int_0^{X(s-)} z\tilde{M}_c(\d s,\d z,\d u) \cr
 \ar\ar\qquad
- \int_0^t X(s-)1_{J_m^c}(s)b(\d s) + \sum_{s\in(0,t]} 1_{J_m}(s)g_s(X(s-)),
 \eeqlb
where
 \beqnn
g_s(x)= \int_0^\infty\int_0^{x} z M_d(\{s\},\d z,\d u) - x\bigg[\Delta b(s) + \int_0^\infty z m_d(\{s\},\d z)\bigg].
 \eeqnn
It is easy to see that $x\mapsto x+g_s(x)$ is an increasing function on $[0,\infty)$. Let $l_s(x_1,x_2)= g_s(x_1) - g_s(x_2)$. Suppose that $\{X_1(t): t\ge0\}$ and $\{X_2(t): t\ge0\}$ are two solutions to \eqref{eq4.6} satisfying $X_1(0)\le X_2(0)$. Let $Y(t)= X_1(t) - X_2(t)$ for $t\ge 0$. Then
 \beqnn
Y(t)\ar=\ar Y(0) + \int_0^t \int_{X_1(s-)\land X_2(s-)}^{X_1(s-)\vee X_2(s-)} (1_{\{Y(s-)>0\}} - 1_{\{Y(s-)<0\}}) W(\d s,\d u) \cr
 \ar\ar
+ \int_0^t\int_0^\infty \int_{X_1(s-)\land X_2(s-)}^{X_1(s-)\vee X_2(s-)} z(1_{\{Y(s-)>0\}} - 1_{\{Y(s-)<0\}}) \tilde{M}_c(\d s,\d z,\d u) \cr
 \ar\ar
- \int_0^tY(s-)1_{J_m^c}(s)b(\d s) + \sum_{s\in (0,t]} l_s(X_1(s-),X_2(s-)) 1_{J_m}(s).
 \eeqnn
By It\^{o}'s formula,
 \beqlb\label{eq4.7}
f_n(Y(t))
 \ar=\ar
\int_0^t Y(s-)f_n''(Y(s-))1_{\{Y(s-)>0\}} c(\d s) + \textrm{local~mart.} \cr
\ar\ar
+ \int_0^t \int_0^\infty\int_{X_1(s-)\land X_2(s-)}^{X_1(s-)\vee X_2(s-)} D_zf_n(Y(s-))1_{\{Y(s-)>0\}} M_c(\d s,\d z,\d u) \cr
 \ar\ar
- \int_0^t f'_n(Y(s-))Y(s-)1_{\{Y(s-)>0\}} b(\d s) \cr
 \ar\ar
+ \sum_{s\in (0,t]}D_{-Y(s-)\Delta b(s)}f_n(Y(s-)) 1_{\{Y(s-)>0\}} 1_{J_b\setminus J_m}(s) \cr
 \ar\ar
+ \sum_{s\in (0,t]} D_{l_s(X_1(s-),X_2(s-))}f_n(Y(s-)) 1_{\{Y(s-)>0\}} 1_{J_m}(s),
 \eeqlb
where $D_zf_n(x)= f_n(x + z) - f_n(x) - f¡ä_n(x)z$. Following the arguments in the proof of Theorem~8.2 in Li (2019+) one sees $|f_n^{\prime\prime}(Y(s-))|\le 2n^{-1}|Y(s-)|^{-1}$ and $|D_zf_n(Y(s-))|\le 2n^{-1}|Y(s-)|^{-1}z^2$. By Lemma~3.1 in Li and Pu (2012) we have, for $s\in J_b\setminus J_m$,
 \beqnn
\big|D_{-Y(s-)\Delta b(s)}f_n(Y(s-))\big|
 \le
2n^{-1}|Y(s-)|^{-1}|Y(s-)\Delta b(s)|^2
 \le
2n^{-1}|Y(s-)||\Delta b(s)|^2
 \eeqnn
and, for $s\in J_m$,
 \beqnn
\ar\ar\big|D_{l_s(X_1(s-),X_2(s-))}f_n(Y(s-))\big| \ccr
 \ar\ar\qquad
\le 2n^{-1}|Y(s-)|^{-1} \bigg(\int_0^\infty \int_{X_2(s-)}^{X_1(s-)} z \tilde{M}_d(\{s\},\d z,\d u) - Y(s-)\Delta b(s)\bigg)^2 \cr
 \ar\ar\qquad
\le 4n^{-1}\bigg[|Y(s-)|^{-1}\bigg(\int_0^\infty \int_{X_2(s-)}^{X_1(s-)} z \tilde{M}_d(\{s\},\d z,\d u)\bigg)^2 + |Y(s-)||\Delta b(s)|^2\bigg]. \qquad\quad
 \eeqnn
It follows that, for $s\in J_b\setminus J_m$,
 \beqnn
\mbf{P}\big[\big|D_{-Y(s-)\Delta b(s)}f_n(Y(s-))\big| 1_{\{Y(s-)>0\}}\big]
 \le
2n^{-1}\mbf{P}(Y(s-)\vee 0)|\Delta b(s)|^2
 \eeqnn
and, for $s\in J_m$,
 \beqnn
\ar\ar\mbf{P}\big[\big|D_{l_s(X_1(s-),X_2(s-))}f_n(Y(s-))\big|1_{\{Y(s-)>0\}}\big] \ccr
 \ar\ar\qquad
\le 4n^{-1}\bigg[\int_0^\infty z^2 m_d(\{s\},\d z) + \mbf{P}(Y(s-)\vee 0)|\Delta b(s)|^2\bigg]. \qquad\quad
 \eeqnn
By taking the expectations in \eqref{eq4.7} we get
 \beqnn
\mbf{P}[f_n(Y(t))]
 \ar\le\ar
\int_0^t\mbf{P}(Y(s-)\vee 0)\|b\|(\d s) + 2n^{-1}c(t) + 2n^{-1}\int_0^t \int_0^\infty z^2 m_c(\d s,\d z) \cr
 \ar\ar
+ \,4n^{-1}\int_0^t\int_0^\infty z^2 m_d(\d s,\d z) + 4n^{-1}\sum_{s\in (0,t]} \mbf{P}(Y(s-)\vee 0)|\Delta b(s)|^2.
 \eeqnn
Then letting $n\to\infty$ gives
 \beqnn
\mbf{P}(Y(t)\vee 0)
 \le
\int_0^t \mbf{P}(Y(s-)\vee 0)\|b\|(\d s).
 \eeqnn
By Gronwall's inequality one can see $\mbf{P}(Y(t)\vee 0)= 0$ for every $t\ge 0$. That proves the desired result. \qed

The following result gives a characterization of the conditional distribution of the jump of the CBVE-process at any moment $t\in J_b\cup J_m$.

\bproposition\label{th4.3} The CBVE-process with transition semigroup $(Q_{r,t})_{t\ge r}$ given by \eqref{eq2.18} and \eqref{eq3.7} has a c\`{a}dl\`{a}g semimartingale realization $\{(X(t),\mcr{F}_t): t\ge 0\}$ with the filtration satisfying the usual hypotheses. For such a realization and $t\in J_b\cup J_m$ we have
 \beqlb\label{eq4.8}
\mbf{P}\big(\e^{-\lambda \Delta X(t)}\big|\mcr{F}_{t-}\big)
 =
\e^{(\lambda-v_{t-,t}(\lambda))X(t-)}, \qquad \lambda\ge 0,
 \eeqlb
where $\Delta X(t)= X(t) - X(t-)$ and
 \beqlb\label{eq4.9}
\lambda-v_{t-,t}(\lambda)
 =
\Delta b(t)\lambda + \int_0^\infty K(\lambda,z) m(\{t\},\d z).
 \eeqlb
\eproposition

\proof Let $\{(X(t), \mcr{G}_t):t\ge 0\}$ be a realization of the CBVE-process defined on a complete probability space $(\Omega, \mcr{F}, \mbf{P})$. In view of \eqref{eq3.7}, for any $\lambda\ge 0$ we have $v_{r,t}(\lambda)\to \lambda$ as $t\downarrow r$. Then \eqref{eq2.18} implies $\lim_{t\downarrow r} Q_{r,t}(x,\d y)= \delta_x(\d y)$ by weak convergence and so $\lim_{t\downarrow r} Q_{r,t}(x,\{y\ge 0: |y-x|\ge \varepsilon\})= 0$ for every $\varepsilon > 0$. By dominated convergence,
 \beqnn
\lim_{t\downarrow r} \mbf{P}\{|X(t)-X(r)|>\varepsilon\}
 =
\lim_{t\downarrow r} \mbf{P}[Q_{r,t}(x, \{y\ge 0: |y-x|\ge \varepsilon\})|_{x=X(r)}] = 0.
 \eeqnn
Then $\{X(t): t\ge 0\}$ is stochastically right continuous. From \eqref{eq3.7} we see $r\mapsto v_{r,t}(\lambda)$ is right-continuous on $[0,t]$, so $\{\e^{-X(r)v_{r,t}(\lambda)}: r\in [0, t]\}$ is stochastically right-continuous. Let $\bar{\mcr{G}}_t$ be the augmentation of $\mcr{G}_t$ and let $\mcr{F}_t= \bar{\mcr{G}}_{t+}$ for $t\ge 0$. The Markov property implies
 \beqlb\label{eq4.10}
\mbf{P}\big[\e^{-\lambda X(t)}\big|\bar{\mcr{G}}_r\big]
 =
\e^{-X(r)v_{r,t}(\lambda)}, \qquad t\ge r\ge 0.
 \eeqlb
This means $\{\e^{-X(r)v_{r,t}(\lambda)}: r\in [0,t]\}$ is a positive bounded martingale, so it has a c\`{a}dl\`{a}g $(\mcr{F}_r)$-martingale modification. By Proposition~\ref{th2.4} we have $v_{r,t}(\lambda)\ge l_{0,t}(\lambda)> 0$ for $\lambda> 0$. Then $\{X(r): r\in [0,t]\}$ has a c\`{a}dl\`{a}g modification. It follows that $\{X(t): t\ge 0\}$ has a c\`{a}dl\`{a}g semimartingale modification; see, e.g., Dellacherie and Meyer (1982, pp.219-221). Using such a modification we can replace $\bar{\mcr{G}}_r$ by $\mcr{F}_r$ in \eqref{eq4.10}. Then $\{(X(t),\mcr{F}_t): t\ge 0\}$ is a c\`{a}dl\`{a}g semimartingale realization of the CBVE-process with the filtration satisfying the usual hypotheses. By letting $r\uparrow t$ in \eqref{eq3.7} and \eqref{eq4.10} we get
 \beqnn
\mbf{P}\big[\e^{-\lambda X(t)}\big|\mcr{F}_{t-}\big]
 =
\e^{-X(t-)v_{t-,t}(\lambda)}, \qquad \lambda\ge 0,
 \eeqnn
where
 \beqnn
v_{t-,t}(\lambda)= (1-\Delta b(t))\lambda - \int_0^\infty K(\lambda,z) m(\{t\},\d z).
 \eeqnn
Then \eqref{eq4.8} follows. \qed

In view of \eqref{eq4.8} and \eqref{eq4.9}, for any $t\in J_b\cup J_m$ it is natural to expect that the jump $\Delta X(t)$ of the CBVE-process should be given by the position at \textit{time} $X(t-)$ of a spectrally positive L\'{e}vy process with L\'{e}vy measure $m(\{t\},\d z)$. It can be realized by an extension of the probability space. For this purpose we first establish a composite L\'{e}vy--It\^{o} representation as follows.

\bproposition\label{th4.4} Let $(\Omega,\mcr{F},\mbf{P})$ be a complete probability space with the sub-$\sigma$-algebra $\mcr{G}\subset \mcr{F}$. Suppose that $(\xi,Z)$ is a random vector taking values in $[0,\infty)\times \mbb{R}$ such that $\xi$ is $\mcr{G}$-measurable and, for every $\lambda\ge 0$,
 \beqlb\label{eq4.11}
\mbf{P}(\e^{-\lambda Z}|\mcr{G})
 =
\exp\bigg\{\xi\bigg[\beta\lambda + \int_0^\infty (\e^{-\lambda z}-1+\lambda z)\gamma(\d z)\bigg]\bigg\},
 \eeqlb
where $\beta\in\mbb{R}$ and $\gamma(\d z)$ is a $\sigma$-finite measure on $(0,\infty)$ satisfying
 \beqnn
\int_0^\infty (z\land z^2) \gamma(\d z)< \infty.
 \eeqnn
Then on an extension of the probability space there exists a Poisson random measure $N(\d z,\d u)$ on $(0,\infty)^2$ with intensity $\gamma(\d z)\d u$ such that $N$ is independent of $\mcr{G}$ and a.s.\
 \beqlb\label{eq4.12}
Z= -\beta\xi + \int_0^\infty\int_0^\xi z \tilde{N}(\d z,\d u).
 \eeqlb
\eproposition

\proof This proof also makes precise the statements of the proposition. Let $\rho(z,u)= (z\land z^2)(1+u^2)^{-1}$ for $z>0$ and $u>0$. Let $M_\rho$ denote the space of all $\sigma$-finite Borel measures $\nu$ on $(0,\infty)^2$ so that
 \beqnn
\int_0^\infty\int_0^\infty \rho(z,u)\nu(\d z,\d u)< \infty.
 \eeqnn
We equip $M_\rho$ with the $\sigma$-algebra $\mcr{M}_\rho$ generated by the mappings $\nu\mapsto \nu((a,\infty)\times B)$ for all $a>0$ and bounded $B\in \mcr{B}(0,\infty)$. It is well-known that there is a spectrally positive L\'{e}vy process $\{Y_s: s\ge 0\}$ such that
 \beqnn
\mbf{E}(\e^{-\lambda Y_s})
 =
\exp\bigg\{s\bigg[\beta\lambda + \int_0^\infty (\e^{-\lambda z}-1+\lambda z)\gamma(\d z)\bigg]\bigg\}, \qquad \lambda\ge 0.
 \eeqnn
By L\'{e}vy--It\^{o} representation, there is a Poisson random measure $G= G(\d z,\d u)$ on $(0,\infty)^2$ with intensity $\gamma(\d z)\d u$ such that
 \beqnn
Y_s= -\beta s + \int_0^\infty\int_0^s z \tilde{G}(\d z,\d u), \qquad s\ge 0.
 \eeqnn
Let $P(s,\d y,\d\nu)$ be the joint distribution of the random vector $(Y_s,G)$ on $\mbb{R}\times M_\rho$. Let $P_1(s,\d y)$ and $P_2(s,\d\nu)$ denote the marginal distributions of $Y_s$ and $G$, respectively. Let $\kappa_1(s,y,\d\nu)$ be a regular conditional distribution of $G$ given $Y_s$. Then $\kappa_1(s,y,\d\nu)$ is a kernel from $[0,\infty)\times \mbb{R}$ to $M_\rho$ and $P(s,\d y,\d\nu)= P_1(s,\d y)\kappa_1(s,y,\d\nu)$. Let $\tilde{\Omega}= \Omega\times M_\rho$ and $\tilde{\mcr{F}}= \mcr{F}\times \mcr{M}_\rho$. Let $\tilde{\mbf{P}}$ be the probability law on $(\tilde{\Omega}, \tilde{\mcr{F}})$ defined by $\tilde{\mbf{P}}(\d\tilde{\omega})= \mbf{P}(\d\omega) \kappa_1(\xi(\omega),Z(\omega),\d\mu)$, where $\tilde{\omega}= (\omega,\mu)\in \tilde{\Omega}$. For any random variable $X$ on $(\Omega, \mcr{F}, \mbf{P})$, write $X(\tilde{\omega})= X(\omega)$ for $\tilde{\omega}= (\omega,\mu)\in \tilde{\Omega}$, which extends $X$ to a random variable on $(\tilde{\Omega}, \tilde{\mcr{F}}, \tilde{\mbf{P}})$. It is easy to see that $\tilde{\mcr{G}}:= \mcr{G}\times \{\emptyset, M_\rho\}\subset \tilde{\mcr{F}}$ and $\xi$ is $\tilde{\mcr{G}}$-measurable as a random variable on $(\tilde{\Omega}, \tilde{\mcr{F}}, \tilde{\mbf{P}})$. Let $N(\tilde{\omega})= \mu$ for $\tilde{\omega}= (\omega,\mu)\in \tilde{\Omega}$. By \eqref{eq4.11} we have $\mbf{P}(Z\in \d y|\mcr{G})= \mbf{P}(Z\in \d y|\xi)= P_1(\xi,\d y)$. From the definition of $\tilde{\mbf{P}}$ it follows that
 \beqnn
\tilde{\mbf{P}}(Z\in \d y, N\in \d\nu|\tilde{\mcr{G}})
 \ar=\ar
\tilde{\mbf{P}}(Z\in \d y, N\in \d\nu|\xi) \cr
 \ar=\ar
P_1(\xi,\d y)\kappa_1(\xi,y,\d\nu)
 =
P(\xi,\d y,\d\nu).
 \eeqnn
Then $N$ is a Poisson random measure on $(0,\infty)^2$ with intensity $\gamma(\d z)\d u$ and \eqref{eq4.12} a.s.\ holds. Let $F$ be a bounded $\mcr{G}$-measurable random variable on $(\Omega, \mcr{F}, \mbf{P})$. For any positive Borel function $f$ on $(0,\infty)^2$ bounded above by $\rho\cdot\const.$, we have
 \beqnn
\ar\ar\tilde{\mbf{P}}\bigg[F\exp\bigg\{-\int_0^\infty\int_0^\infty f(z,u) N(\d z,\d u)\bigg\}\bigg] \cr
 \ar\ar\qquad
= \tilde{\mbf{P}}\bigg[F\tilde{\mbf{P}}\bigg(\exp\bigg\{-\int_0^\infty\int_0^\infty f(z,u) N(\d z,\d u)\bigg\}\bigg|\tilde{\mcr{G}}\bigg)\bigg] \cr
 \ar\ar\qquad
= \tilde{\mbf{P}}\bigg[F\int_{M_\rho}\exp\bigg\{-\int_0^\infty\int_0^\infty f(z,u) \nu(\d z,\d u)\bigg\}P_2(\xi,\d\nu)\bigg] \cr
 \ar\ar\qquad
= \tilde{\mbf{P}}(F)\exp\bigg\{-\int_0^\infty \gamma(\d z)\int_0^\infty (1-\e^{-f(z,u)})\d u\bigg\} \cr
 \ar\ar\qquad
= \tilde{\mbf{P}}(F) \tilde{\mbf{P}}\bigg[\exp\bigg\{-\int_0^\infty\int_0^\infty f(z,u) N(\d z,\d u)\bigg\}\bigg].
 \eeqnn
Then $N$ is independent of $\tilde{\mcr{G}}$ on $(\tilde{\Omega}, \tilde{\mcr{F}}, \tilde{\mbf{P}})$. \qed


\btheorem\label{th4.5} There is a pathwise unique solution $\{X(t): t\ge 0\}$ to \eqref{eq4.1} and the solution is a CBVE-process with transition semigroup $(Q_{r,t})_{t\ge r}$ defined by \eqref{eq2.18} and \eqref{eq3.7}.
\etheorem

\proof \textit{Step~1.} Consider the case where \eqref{eq4.4} holds and $b(t)= 0$ for every $t\ge 0$. Let $(v_{r,t})_{t\ge r}$ be the conservative cumulant semigroup defined by \eqref{eq3.7} in this special case. Suppose that $\{(X(t),\mcr{F}_t): t\ge 0\}$ is a the realization of the corresponding CBVE-process provided by Proposition~\ref{th4.3}. Then the process is actually a martingale by Proposition~\ref{th3.5}. Let $N_0(\d s,\d z)$ be the optional random measure on $(0,\infty)\times \mbb{R}$ by
 \beqnn
N_0(\d s,\d z):= \sum_{s>0}1_{\{\Delta X(s)\neq0\}} \delta_{(s,\Delta X(s))}(\d s,\d z).
 \eeqnn
Let $\hat{N}_0(\d s,\d z)$ denote the predictable compensator of $N_0(\d s,\d z)$ and let $\tilde{N}_0(\d s,\d z)= N_0(\d s,\d z) - \hat{N}_0(\d s,\d z)$ be the compensated measure. We can write
 \beqlb\label{eq4.13}
X(t)= X(0) + M(t) + \int_0^t \int_{\mbb{R}} z \tilde{N}_0(\d s,\d z),
 \eeqlb
where $\{M(t): t\ge 0\}$ is a continuous local martingale. Let $\{C(t): t\ge 0\}$ be its quadratic variation process. Let $f(x,\lambda)=\e^{-x\lambda}$ for $x,\lambda\ge 0$. Then
 \beqnn
f'_1(x,\lambda)= -\lambda f(x,\lambda), ~ f'_2(x,\lambda)= -xf(x,\lambda), ~ f''_{11}(x,\lambda)= \lambda^2 f(x,\lambda).
 \eeqnn
By It\^{o}'s formula, for $t\ge r\ge 0$ and $\lambda\ge 0$,
 \beqlb\label{eq4.14}
\e^{-X(t)\lambda}\ar=\ar \e^{-X(r)v_{r,t}(\lambda)} + \int_r^t f_1'(X(s-), v_{s-,t}(\lambda))\d X(s) + \int_r^t f_2'(X(s-), v_{s-,t}(\lambda))\d v_{s,t}(\lambda) \cr
 \ar\ar
+\, \frac{1}{2}\int_r^t f_{11}''(X(s-), v_{s-,t}(\lambda)) \d C(s) \cr
 \ar\ar
+ \sum_{s\in (r,t]\cap J_m}\Big[f(X(s),v_{s,t}(\lambda)) - f(X(s-),v_{s-,t}(\lambda)) \cr
 \ar\ar\qqquad\quad
-\, f_1'(X(s-),v_{s-,t}(\lambda))\Delta X(s) - f_2'(X(s-),v_{s-,t}(\lambda))\Delta v_{s,t}(\lambda)\Big] \ccr
 \ar\ar
+ \sum_{s\in (r,t]\setminus J_m}\Big[f(X(s), v_{s,t}(\lambda))-f(X(s-), v_{s,t}(\lambda)) -f_1'(X(s-), v_{s,t}(\lambda))\Delta X(s)\Big] \cr
 \ar=\ar
\e^{-X(r)v_{r,t}(\lambda)} - \int_r^t \e^{-X(s-)v_{s-,t}(\lambda)}v_{s-,t}(\lambda)\d X(s) \cr
 \ar\ar
- \int_r^t \e^{-X(s-)v_{s-,t}(\lambda)}X(s-)v_{s,t}(\lambda)^2 c(\d s) \cr
 \ar\ar
- \int_r^t\int_0^\infty \e^{-X(s-)v_{s-,t}(\lambda)}X(s-) K(v_{s,t}(\lambda),z) m(\d s,\d z) \cr
 \ar\ar
+\, \frac{1}{2} \int_r^t \e^{-X(s-)v_{s-,t}(\lambda)}v_{s,t}(\lambda)^2C(\d s) + \sum_{s\in (r,t]\cap J_m} \Big[\e^{-X(s)v_{s,t}(\lambda)} - \e^{-X(s-)v_{s-,t}(\lambda)} \cr
 \ar\ar\qqquad
+\, \e^{-X(s-)v_{s-,t}(\lambda)}v_{s-,t}(\lambda)\Delta X(s) + \e^{-X(s-)v_{s-,t}(\lambda)}X(s-) \Delta v_{s,t}(\lambda)\Big] \cr
 \ar\ar
+ \sum_{s\in(r,t]\setminus J_m}\Big[\e^{-X(s)v_{s,t}(\lambda)} - \e^{-X(s-)v_{s,t}(\lambda)} + \e^{-X(s-)v_{s,t}(\lambda)}v_{s,t}(\lambda)\Delta X(s)\Big] \cr
 \ar=\ar
\e^{-X(r)v_{r,t}(\lambda)} + \int_0^r \e^{-X(s-)v_{s-,t}(\lambda)}v_{s-,t}(\lambda)\d X(s) \cr
 \ar\ar
+ \int_0^r \e^{-X(s-)v_{s-,t}(\lambda)}X(s-)v_{s,t}(\lambda)^2 c(\d s) \cr
 \ar\ar
+ \int_0^r\int_0^\infty \e^{-X(s-)v_{s-,t}(\lambda)}X(s-) K(v_{s,t}(\lambda),z) m(\d s,\d z) + Z(t) \cr
 \ar\ar
-\, \frac{1}{2} \int_0^r \e^{-X(s-)v_{s-,t}(\lambda)}v_{s,t}(\lambda)^2C(\d s) - \sum_{s\in (0,r]\cap J_m} \Big[\e^{-X(s)v_{s,t}(\lambda)} - \e^{-X(s-)v_{s-,t}(\lambda)} \cr
 \ar\ar\qqquad
+\, \e^{-X(s-)v_{s-,t}(\lambda)}v_{s-,t}(\lambda)\Delta X(s) + \e^{-X(s-)v_{s-,t}(\lambda)}X(s-)\Delta v_{s,t}(\lambda)\Big] \cr
 \ar\ar
- \sum_{s\in(0,r]\setminus J_m}\e^{-X(s-)v_{s,t}(\lambda)}\Big[\e^{-\Delta X(s)v_{s,t}(\lambda)} - 1 + v_{s,t}(\lambda)\Delta X(s)\Big],
 \eeqlb
where
 \beqnn
Z(t)\ar=\ar -\int_0^t \e^{-X(s-)v_{s-,t}(\lambda)}v_{s-,t}(\lambda)\d X(s) - \int_0^t \e^{-X(s-)v_{s-,t}(\lambda)}X(s-)v_{s,t}(\lambda)^2 c(\d s) \cr
 \ar\ar
- \int_0^t\int_0^\infty \e^{-X(s-)v_{s-,t}(\lambda)}X(s-) K(v_{s,t}(\lambda),z) m(\d s,\d z) \cr
 \ar\ar
+ \frac{1}{2} \int_0^t \e^{-X(s-)v_{s-,t}(\lambda)}v_{s,t}(\lambda)^2C(\d s) + \sum_{s\in(0,t]\cap J_m} \Big[\e^{-X(s)v_{s,t}(\lambda)} - \e^{-X(s-)v_{s-,t}(\lambda)} \cr
 \ar\ar\qqquad
+\, \e^{-X(s-)v_{s-,t}(\lambda)} v_{s-,t}(\lambda)\Delta X(s) + \e^{-X(s-)v_{s-,t}(\lambda)}X(s-)\Delta v_{s,t}(\lambda)\Big] \cr
 \ar\ar
+ \sum_{s\in(0,t]\setminus J_m}\e^{-X(s-)v_{s,t}(\lambda)}\Big[\e^{-\Delta X(s)v_{s,t}(\lambda)} - 1 + v_{s,t}(\lambda)\Delta X(s)\Big],
 \eeqnn
Taking the conditional expectation in \eqref{eq4.14} we obtain
 \beqnn
\mbf{P}[-Z(t)|\mcr{F}_r]
 \ar=\ar
\int_0^r \e^{-X(s-)v_{s-,t}(\lambda)}X(s-)v_{s,t}(\lambda)^2 c(\d s) + \mbox{mart.} \cr
 \ar\ar
+ \int_0^r\int_0^\infty \e^{-X(s-)v_{s-,t}(\lambda)}X(s-) K(v_{s,t}(\lambda),z) m(\d s,\d z) \cr
 \ar\ar
-\, \frac{1}{2} \int_0^r \e^{-X(s-)v_{s-,t}(\lambda)}v_{s,t}(\lambda)^2C(\d s) \cr
 \ar\ar
- \sum_{s\in (0,r]\cap J_m} \e^{-X(s-)v_{s-,t}(\lambda)}\big[\Delta X(s)v_{s-,t}(\lambda) + X(s-)\Delta v_{s,t}(\lambda)\big] \cr
 \ar\ar
- \sum_{s\in(0,r]\setminus J_m}\e^{-X(s-)v_{s,t}(\lambda)}\Big[\e^{-\Delta X(s)v_{s,t}(\lambda)} - 1 + v_{s,t}(\lambda)\Delta X(s)\Big] \cr
 \ar=\ar
\int_0^r \e^{-X(s-)v_{s-,t}(\lambda)}X(s-)v_{s,t}(\lambda)^2 c(\d s) + \mbox{mart.} \cr
 \ar\ar
+ \int_0^r\int_0^\infty \e^{-X(s-)v_{s-,t}(\lambda)}X(s-) K(v_{s,t}(\lambda),z) m_c(\d s,\d z) \cr
 \ar\ar
-\, \frac{1}{2} \int_0^r \e^{-X(s-)v_{s-,t}(\lambda)}v_{s,t}(\lambda)^2C(\d s) \cr
 \ar\ar
- \sum_{s\in (0,r]\cap J_m} \e^{-X(s-)v_{s-,t}(\lambda)}\Delta X(s)v_{s-,t}(\lambda) \cr
 \ar\ar
- \sum_{s\in(0,r]\setminus J_m}\e^{-X(s-)v_{s,t}(\lambda)}\Big[\e^{-\Delta X(s)v_{s,t}(\lambda)} - 1 + v_{s,t}(\lambda)\Delta X(s)\Big],
 \eeqnn
Then the uniqueness of canonical decompositions of martingales yields
 \beqnn
\d C(s)= 2X(s-)c(\d s), \quad 1_{J_m^c}(s)\hat{N_0}(\d s,\d z)= X(s-)m_c(\d s,\d z).
 \eeqnn
By El~Karoui and M\'{e}l\'{e}ard (1990, Theorem~III.6), on an extension of the original probability space there exists a Gaussian white noise $W(\d s,\d u)$ on $(0,\infty)^2$ with intensity $2c(\d s)\d u$ such that
 \beqnn
M(t)= \int_0^t\int_0^{X(s-)}W(\d s,\d u).
 \eeqnn
By Kabanov et al.\ (1981, Theorem~1), on a further extension of the original probability space we can define a Poisson random measure $M_c(\d s,\d z,\d u)$ with intensity $m_c(\d s,\d z)\d u$ so that
 \beqnn
\int_0^t\int_0^\infty z1_{J_m^c}(s)\tilde{N}_0(\d s,\d z)
 =
\int_0^t\int_0^\infty\int_0^{X(s-)} z\tilde{M}_c(\d s,\d z,\d u);
 \eeqnn
see also El~Karoui and Lepeltier (1977). By \eqref{eq4.13} we see the process a.s.\ makes a jump at time $s\in J_m$ with the representation:
 \beqnn
\Delta X(s)= \int_{\mbb{R}} z \tilde{N}_0(\{s\},\d z).
 \eeqnn
From Proposition~\ref{th4.3} it follows that
 \beqnn
\mbf{P}\big(\e^{-\lambda \Delta X(s)}\big|\mcr{F}_{s-}\big)
 =
\e^{(\lambda-v_{s-,s}(\lambda))X(s-)}, \qquad \lambda\ge 0,
 \eeqnn
where
 \beqnn
\lambda-v_{s-,s}(\lambda)
 =
\int_0^\infty (\e^{-\lambda z}-1+\lambda z) m_d(\{s\},\d z).
 \eeqnn
By Proposition~\ref{th4.4} we can make another extension of the probability space and define a Poisson random measure $N_s(\d z,\d u)$ on $(0,\infty)^2$ with intensity $m_d(\{s\},\d z)\d u$ such that $N_s$ is independent of $\mcr{F}_{s-}$ and
 \beqnn
\Delta X(s)= \int_0^\infty\int_0^{X(s-)} z \tilde{N}_s(\d z,\d u).
 \eeqnn
Let $M_d(\d s,\d z,\d u)$ be the random measure on $(0,\infty)^3$ defined by
 \beqnn
M_d((0,t]\times A)
 =
\sum_{s\in(0,t]\cap J_m} N_s(A), \qquad t\ge 0, ~ A\in \mcr{B}((0,\infty)^2).
 \eeqnn
Then $M_d(\d s,\d z,\d u)$ is a Poisson random measure with intensity $m_d(\d s,\d z)\d u$. From \eqref{eq4.13} we see $\{X(t): t\ge 0\}$ is a solution to \eqref{eq4.1} with $b=0$.

\textit{Step~2.} Let us show that the noises $W$, $M_c$ and $M_d$ constructed in the first step are independent. Let $g$ be a positive continuous function on $(0,\infty)$ and let $f$ and $h$ be positive continuous functions on $(0,\infty)^2$. We assume all those functions have compact supports. For $t\ge 0$ write
 \beqnn
F_t(g,f,h)\ar=\ar \int_0^t\int_0^\infty g(u) W(\d s,\d u) + \int_0^t\int_0^\infty\int_0^\infty f(z,u) M_c(\d s,\d z,\d u) \cr
 \ar\ar\qqquad
+ \int_0^t\int_0^\infty\int_0^\infty h(z,u) M_d(\d s,\d z,\d u).
 \eeqnn
By It\^{o}'s formula, we have
 \beqnn
\e^{-F_t(g,f,h)}\ar=\ar 1 + \int_0^t \e^{-F_{s-}(g,f,h)}c(\d s) \int_0^\infty g(u)^2\d u + \mbox{mart.} \cr
 \ar\ar
+ \int_0^t\int_0^\infty\int_0^\infty \big[\e^{-F_{s-}(g,f,h)-f(z,u)} - \e^{-F_{s-}(g,f,h)}\big] M_c(\d s,\d z,\d u) \cr
 \ar\ar
+ \sum_{s\in(0,t]\cap J_m} \big[\e^{-F_s(g,f,h)} - \e^{-F_{s-}(g,f,h)}\big] \cr
 \ar=\ar
1 + \int_0^t \e^{-F_{s-}(g,f,h)}c(\d s) \int_0^\infty g(u)^2\d u + \mbox{mart.} \cr
 \ar\ar
+ \int_0^t\int_0^\infty\int_0^\infty \e^{-F_{s-}(g,f,h)}(\e^{-f(z,u)}-1) M_c(\d s,\d z,\d u) \cr
 \ar\ar
+ \sum_{s\in(0,t]\cap J_m} \e^{-F_{s-}(g,f,h)}(\e^{-M_d(s,h)}-1) \cr
 \ar=\ar
1 + \int_0^t \e^{-F_{s-}(g,f,h)} c(\d s) \int_0^\infty g(u)^2\d u + \mbox{mart.} \cr
 \ar\ar
+ \int_0^t\int_0^\infty \e^{-F_{s-}(g,f,h)} m_c(\d s,\d z)\int_0^\infty(\e^{-f(z,u)}-1)\d u \cr
 \ar\ar
+ \sum_{s\in(0,t]\cap J_m} \e^{-F_{s-}(g,f,h)}(\e^{-m_d(s,h)}-1),
 \eeqnn
where
 \beqnn
M_d(s,h)= \int_0^\infty\int_0^\infty h(z,u) M_d(\{s\},\d z,\d u)
 \eeqnn
and
 \beqnn
m_d(s,h)= \int_0^\infty m_d(\{s\},\d z)\int_0^\infty (1-\e^{-h(z,u)})\d u.
 \eeqnn
Then we can use Proposition~\ref{th2.2} to obtain
 \beqnn
\mbf{P}\big[\e^{-F_t(g,f,h)}\big]
 \ar=\ar
\exp\bigg\{\int_0^t c(\d s) \int_0^\infty g(u)^2\d u - \sum_{s\in(0,t]\cap J_m} m_d(s,h) \cr
 \ar\ar\qquad
- \int_0^t\int_0^\infty m_c(\d s,\d z)\int_0^\infty(1-\e^{-f(z,u)})\d u\bigg\}.
 \eeqnn
That gives the independence of $W$, $M_c$ and $M_d$.

\textit{Step~3.} Consider the case where \eqref{eq4.4} holds. Let $(v_{r,t})_{t\ge r}$ be the conservative cumulant semigroup defined by \eqref{eq3.7}. Let $\zeta$ be the c\`{a}dl\`{a}g function on $[0,\infty)$ such that $\zeta_c(t)= -b_c(t)$ and $\Delta\zeta(t)= \log[1-\Delta b(t)]$ for every $t\ge 0$. By Theorem~\ref{th3.6} we can define another cumulant semigroup $(u_{r,t})_{t\ge r}$ by \eqref{eq3.16} and $r\mapsto u_{r,t}(\lambda)$ is the unique bounded positive solution to
 \beqnn
u_{r,t}(\lambda)= \lambda - \int_r^tu_{s,t}(\lambda)^2 \e^{-\zeta(s)}c(\d s) - \int_r^t\int_0^\infty K(u_{s,t}(\lambda),z) \e^{\zeta(s-)}m(\d s,\e^{\zeta(s)}\d z).
 \eeqnn
Let $W(\d s,\d u)$ and $M(\d s,\d z,\d u)$ be given as in the introduction. One can see that $W_0(\d s,\d u):= \e^{-\zeta(s)}W(\d s,\e^{\zeta(s-)}\d u)$ is a Gaussian white noise on $(0,\infty)^2$ with intensity $2\e^{-\zeta(s)}c(\d s)\d u$ and $M_0(\d s,\d z,\d u):= M(\d s,\e^{\zeta(s)}\d z,\e^{\zeta(s-)}\d u)$ is a Poisson random measure on $(0,\infty)^3$ with intensity $\e^{\zeta(s-)}m(\d s,\e^{\zeta(s)}\d z)\d u$. By Proposition~\ref{th4.2} and the first step of the proof, we can construct a CBVE-process with cumulant semigroup $(u_{r,t})_{t\ge r}$ by the pathwise unique solution to
 \beqnn
Z(t)= X(0) + \int_0^t\int_0^{Z(s-)}W_0(\d s,\d u) + \int_0^t\int_0^\infty \int_0^{Z(s-)} z\tilde{M}_0(\d s,\d z,\d u).
 \eeqnn
It is easy to see that $t\mapsto X(t):= \e^{\zeta(t)}Z(t)$ is a CBVE-process with cumulant semigroup $(v_{r,t})_{t\ge r}$. By integration by parts we have
 \beqnn
X(t)\ar=\ar X(0) + \int_0^t Z(s-)\d\e^{\zeta(s)} + \int_0^t\e^{\zeta(s)}\d Z(s) \cr
 \ar=\ar
X(0) - \int_0^t\e^{\zeta(s-)}Z(s-) b(\d s) + \int_0^t\int_0^{\e^{-\zeta(s-)}X(s-)} \e^{\zeta(s)} W_0(\d s,\d u) \cr
 \ar\ar
+ \int_0^t\int_0^\infty\int_0^{\e^{-\zeta(s-)}X(s-)} \e^{\zeta(s)} z\tilde{M}_0(\d s,\d z,\d u) \cr
 \ar=\ar
X(0) - \int_0^t X(s-) b(\d s) + \int_0^t\int_0^{X(s-)} \e^{\zeta(s)} W_0(\d s,\e^{-\zeta(s-)}\d u) \cr
 \ar\ar
+ \int_0^t\int_0^\infty\int_0^{X(s-)} z\tilde{M}_0(\d s,\e^{-\zeta(s)}\d z,\e^{-\zeta(s-)}\d u).
 \eeqnn
Then $\{X(t): t\ge 0\}$ solves \eqref{eq4.1}.

\textit{Step~4.} Let us consider the general case. By Proposition~\ref{th4.2} and the second step of the proof, for each $k\ge 1$ we can construct a CBVE-process $\{X_k(t): t\ge 0\}$ by the pathwise unique solution to:
 \beqlb\label{eq4.15}
X(t)\ar=\ar X(0) + \int_0^t\int_0^{X(s-)}W(\d s,\d u) + \int_0^t\int_0^\infty \int_0^{X(s-)} (z\land k) \tilde{M}(\d s,\d z,\d u) \cr
 \ar\ar\qquad
- \int_0^t X(s-) b(\d s) - \int_0^t\int_k^\infty X(s-) (z-k) m(\d s,\d z).
 \eeqlb
The cumulant semigroup $(v_{r,t}^{(k)})_{t\ge r}$ of $\{X_k(t): t\ge 0\}$ is defined by:
 \beqlb\label{eq4.16}
v_{r,t}(\lambda)\ar=\ar \lambda - \int_r^t v_{s,t}(\lambda)^2c(\d s) - \int_r^t\int_0^\infty K(v_{s,t}(\lambda),z\land k) m(\d s,\d z) \cr
 \ar\ar
- \int_r^t v_{s,t}(\lambda)b(\d s) - \int_r^t\int_k^\infty v_{s,t}(\lambda) (z-k) m(\d s,\d z).
 \eeqlb
We can rewrite \eqref{eq4.15} into the equivalent form:
 \beqlb\label{eq4.17}
X(t)\ar=\ar X(0) + \int_0^t\int_0^{X(s-)}W(\d s,\d u) + \int_0^t\int_0^k \int_0^{X(s-)} z \tilde{M}(\d s,\d z,\d u) \cr
 \ar\ar\qquad
- \int_0^t X(s-) b(\d s) - \int_0^t\int_k^\infty X(s-) z m(\d s,\d z) \cr
 \ar\ar\qquad
+\, k\int_0^t\int_k^\infty \int_0^{X(s-)} M(\d s,\d z,\d u).
 \eeqlb
This is obtained from \eqref{eq4.1} by modifying the magnitudes of the large jumps. Let $\zeta_{0,k}=0$ and for $i\ge 0$ inductively define
 \beqnn
\zeta_{i+1,k}= \inf\bigg\{t> 0: \int_{\zeta_{i,k}}^{\zeta_{i,k}+t} \int_k^\infty \int_0^{X_k(s-)} M(\d s,\d z,\d u)\ge 1\bigg\}.
 \eeqnn
From \eqref{eq4.17} it is easy to see that $X_{k+1}(t)= X_k(t)$ for $0\le t<\zeta_{1,k}$ and $X_{k+1}(\zeta_{1,k})\ge X_k(\zeta_{1,k})$. By applying Proposition~\ref{th4.2} successively at the stopping times $\eta_{n,k}:= \sum_{i=1}^n \zeta_{i,k}$, $n\ge 1$ we infer that $X_{k+1}(t)\ge X_k(t)$ for every $t\ge 0$. For $r,x\ge 0$ let $\{X_k(r,x,t): t\ge r\}$ be the pathwise unique solution to:
 \beqnn
X(t)\ar=\ar x + \int_r^t\int_0^{X(s-)}W(\d s,\d u) + \int_r^t\int_0^\infty \int_0^{X(s-)} (z\land k) \tilde{M}(\d s,\d z,\d u) \cr
 \ar\ar\qquad
- \int_r^t X(s-) b(\d s) - \int_r^t\int_k^\infty X(s-) (z-k) m(\d s,\d z).
 \eeqnn
By the preceding arguments we have $X_{k+1}(r,x,t)\ge X_k(r,x,t)$, which implies
 \beqnn
v_{r,t}^{(k+1)}(\lambda)\ar=\ar -\log\mbf{P}\exp\{-\lambda X_{k+1}(r,1,t)\} \cr
 \ar\ge\ar
-\log\mbf{P}\exp\{-\lambda X_k(r,1,t)\}= v_{r,t}^{(k)}(\lambda).
 \eeqnn
From \eqref{eq4.16} we see that $r\mapsto v_{r,t}(\lambda):= \uparrow \lim_{k\uparrow \infty} v_{r,t}^{(k)}(\lambda)$ is the unique bounded positive solution to \eqref{eq3.7}.
Observe that $\zeta_{1,k}\ge \tau_{k/2}:= \inf\{t\ge 0: X_k(t)\ge k/2\}$. By Proposition~\ref{th4.1} it is easy to show that $\lim_{k\to \infty}\tau_{k/2}= \infty$. Then $\{X_k(t): t\ge 0\}$ converges increasingly as $k\to \infty$ to a c\`{a}dl\`{a}g process $\{X(t): t\ge 0\}$ which is a CBVE-process with cumulant semigroup $(v_{r,t})_{t\ge r}$. From \eqref{eq4.15} we infer that $\{X(t): t\ge 0\}$ is a solution to \eqref{eq4.1}. The pathwise uniqueness of the solution holds by Proposition~\ref{th4.2}. \qed


\section{Extensions to the general case}

 \setcounter{equation}{0}

In this section, we extend the results established in the last two sections to the general equations \eqref{eq1.3} and \eqref{eq1.5}. Suppose that $(b_1,c,m)$ are admissible parameters defined as in the introduction.

\bproposition\label{th5.1} For any $\lambda> 0$ and $t\ge 0$ the uniqueness of bounded positive solutions holds for \eqref{eq1.3}. \eproposition

\proof Suppose that both $r\mapsto v_{r,t}(\lambda)$ and $r\mapsto w_{r,t}(\lambda)$ are bounded positive solutions to \eqref{eq1.3}. Then $v_{r,t}(\lambda)\land w_{r,t}(\lambda)\ge l_{0,t}(\lambda)> 0$ by Corollary~\ref{th2.3} and Proposition~\ref{th2.4}. From \eqref{eq1.3} it follows that
 \beqnn
|v_{r,t}(\lambda) - w_{r,t}(\lambda)|
 \ar\le\ar
\int_r^t |v_{s,t}(\lambda) - w_{s,t}(\lambda)|\|b_1\|(\d s) + \int_r^t |v_{s,t}(\lambda)^2 - w_{s,t}(\lambda)^2| c(\d s) \cr
 \ar\ar
+ \int_r^t\int_0^1 |K(v_{s,t}(\lambda),z) - K(w_{s,t}(\lambda),z)| m(\d s,\d z) \cr
 \ar\ar
+ \int_r^t\int_1^\infty |\e^{-v_{s,t}(\lambda)z} - \e^{-w_{s,t}(\lambda)z}| m(\d s,\d z) \cr
 \ar\le\ar
\int_r^t |v_{s,t}(\lambda) - w_{s,t}(\lambda)|\|b_1\|(\d s) + 2U_{0,t}(\lambda) \int_r^t |v_{s,t}(\lambda) - w_{s,t}(\lambda)| c(\d s) \cr
 \ar\ar
+\, 2U_{0,t}(\lambda) \int_r^t\int_0^1|v_{s,t}(\lambda) - w_{s,t}(\lambda)|z^2m(\d s,\d z) \cr
 \ar\ar
+ \int_r^t\int_1^\infty|v_{s,t}(\lambda) - w_{s,t}(\lambda)| z\e^{-l_{0,t}(\lambda)z}m(\d s,\d z).
 \eeqnn
By Gronwall's inequality we have $|v_{r,t}(\lambda) - w_{r,t}(\lambda)|= 0$ for $r\in [0,t]$. \qed

\smallskip\noindent\textit{Proof of Theorems~\ref{th1.1} and~\ref{th1.3}.~} By Theorem~\ref{th3.3}, for each $k\ge 1$ there is a cumulant semigroup $(v_{r,t}^{(k)})_{t\ge r}$ defined by the integral evolution equation:
 \beqlb\label{eq5.1}
v_{r,t}(\lambda)\ar=\ar \lambda - \int_r^t v_{s,t}(\lambda)^2c(\d s) - \int_r^t\int_0^\infty K(v_{s,t}(\lambda),z\land k) m(\d s,\d z) \cr
 \ar\ar
- \int_r^t v_{s,t}(\lambda) b_1(\d s) + \int_r^t\int_1^\infty v_{s,t}(\lambda)(z\land k) m(\d s,\d z).
 \eeqlb
By Theorem~\ref{th4.5} we can construct a CBVE-process $\{X_k(t): t\ge 0\}$ with cumulant semigroup $(v_{r,t}^{(k)})_{t\ge r}$ by the pathwise unique solution to:
 \beqnn
X(t)\ar=\ar X(0) + \int_0^t\int_0^{X(s-)}W(\d s,\d u) + \int_0^t\int_0^\infty\int_0^{X(s-)} (z\land k) \tilde{M}(\d s,\d z,\d u) \cr
 \ar\ar
- \int_0^tX(s-)b_1(\d s) + \int_0^t\int_1^\infty X(s-)(z\land k) m(\d s,\d z).
 \eeqnn
The above stochastic equation is equivalent to
 \beqlb\label{eq5.2}
X(t)\ar=\ar X(0) + \int_0^t\int_0^{X(s-)}W(\d s,\d u) + \int_0^t\int_0^1\int_0^{X(s-)} z \tilde{M}(\d s,\d z,\d u) \cr
 \ar\ar
- \int_0^tX(s-)b_1(\d s) + \int_0^t\int_1^\infty\int_0^{X(s-)} (z\land k) M(\d s,\d z,\d u).
 \eeqlb
Let $\zeta_{1,k}$ and $\tau_{k/2}$ be defined as in the last step of the proof of Theorem~\ref{th4.5}. By the arguments in that proof we have $X_{k+1}(t)= X_k(t)$ for $0\le t<\zeta_{1,k}$ and both $\{X_k(t): t\ge 0\}$ and $(v_{r,t}^{(k)})_{t\ge r}$ are increasing in $k\ge 1$. By Proposition~\ref{th2.4} we have $l_{0,t}(\lambda)\le v_{r,t}^{(k)}(\lambda)\le U_{0,t}(\lambda)$. Then for $\lambda>0$ the limit $v_{r,t}(\lambda):= \uparrow \lim_{\uparrow\to \infty} v_{r,t}^{(k)}(\lambda)$ exists and strictly positive. By letting $k\to \infty$ in \eqref{eq5.1} we see $r\mapsto v_{r,t}(\lambda)$ is a solution to \eqref{eq1.3}. The uniqueness of the solution is guaranteed by Proposition~\ref{th5.1}. Clearly, the family $(v_{r,t})_{t\ge r}$ is a cumulant semigroup. It is easy to see that $\lim_{k\to \infty}\zeta_{1,k}= \tau_\infty:= \lim_{k\to \infty}\tau_{k/2}$. Let $\{X(t): t\ge 0\}$ be the c\`{a}dl\`{a}g process such that $X(t)= X_k(t)$ for $0\le t<\zeta_{1,k}$ and $X(t)= \infty$ for $t\ge \tau_\infty$. Then $\{X(t): t\ge 0\}$ is a CBVE-process with cumulant semigroup $(v_{r,t})_{t\ge r}$. From \eqref{eq5.2} we see that $\{X(t): t\ge 0\}$ is a solution to \eqref{eq1.5}. The pathwise uniqueness for \eqref{eq1.5} follows from that for \eqref{eq5.2}.

\smallskip\noindent\textit{Proof of Theorem~\ref{th1.2}.~} It is easy to see that $r\mapsto v_{r,t}(0)$ is indeed a bounded positive solution to \eqref{eq1.3} with $\lambda= 0$. Suppose that $r\mapsto u_{r,t}(0)$ is another positive solution to \eqref{eq1.3} with $\lambda= 0$ and $u_{r,t}(0)> 0$ for some $r\in [0,t]$. Let $t_0= \inf\{r\in [0,t]: u_{r,t}(0)= 0\}$. We clearly have $u_{r,t}(0)= 0$ for $r\in [t_0,t]$, and hence $u_{t_0-,t}(0)= 0$ by \eqref{eq1.3}. Then for any $\lambda> 0$ we can choose $r_0\in [0,t_0)$ so that $u_{r,t}(0)\le l_{0,t}(\lambda)\le v_{r,t}(\lambda)$ when $r\in [r_0,t_0)$. The definition of $t_0$ yields the existence of some $t_1\in [r_0,t_0)$ so that $0< u_{t_1,t}(0)\le v_{t_1,t}(\lambda)$. For $r\in [0,t_1]$ we see from \eqref{eq1.3} that
 \beqnn
u_{r,t}(0)\ar=\ar u_{t_1,t}(0) - \int_r^{t_1} u_{s,t}(0) b_1(\d s) - \int_r^{t_1} u_{s,t}(0)^2 c(\d s) \cr
 \ar\ar\quad
- \int_r^{t_1}\int_0^\infty K_1(u_{s,t}(0),z) m(\d s,\d z).
 \eeqnn
By the uniqueness of the solution we have $u_{r,t}(0)= v_{r,t_1}(u_{t_1,t}(0))\le v_{r,t_1}(v_{t_1,t}(\lambda))= v_{r,t}(\lambda)$. Then $u_{r,t}(0)\le v_{r,t}(\lambda)$ for every $r\in [0,t]$, implying $u_{r,t}(0)\le \lim_{\lambda\downarrow 0} v_{r,t}(\lambda)= v_{r,t}(0)$ for every $r\in [0,t]$. Similarly, one can show $r\mapsto v_{r,t}(\infty)$ is the smallest positive solution to \eqref{eq1.3} with $\lambda= \infty$. \qed

\end{document}